\documentclass[preprint,12pt,fleqn]{elsarticle}

\usepackage{amsmath,amssymb}
\usepackage{cases}

\usepackage{url}
\biboptions{sort&compress}
\usepackage{graphicx}
\usepackage{color}
\usepackage{ulem}
\usepackage{multirow}

\begin{document}

\begin{frontmatter}



\title{Higher-continuity s-version of finite element method with B-spline functions}


\author[aff:1]{Nozomi Magome}
\ead{magome.nozomi.sw@alumni.tsukuba.ac.jp}
\author[aff:2]{Naoki Morita}
\ead{nmorita@kz.tsukuba.ac.jp}
\author[aff:3]{Shigeki Kaneko}
\ead{s_kaneko@save.sys.t.u-tokyo.ac.jp}
\author[aff:2]{Naoto Mitsume\corref{cor1}}
\ead{mitsume@kz.tsukuba.ac.jp}

\cortext[cor1]{Corresponding author.}

\affiliation[aff:1]{organization={University of Tsukuba, Degree Programs in Systems and Information Engineering},
            addressline={Tennodai 1-1-1}, 
            city={Tsukuba},
            postcode={3058573}, 
            state={Ibaraki},
            country={Japan}}

\affiliation[aff:2]{organization={University of Tsukuba, Institute of Systems and Information Engineering},
            addressline={Tennodai 1-1-1}, 
            city={Tsukuba},
            postcode={3058573}, 
            state={Ibaraki},
            country={Japan}}

\affiliation[aff:3]{organization={The University of Tokyo, Department of Systems Innovation, School of Engineering},
            addressline={Hongo 7-3-1}, 
            city={Bunkyo-ku},
            postcode={1138656}, 
            state={Tokyo},
            country={Japan}}

\begin{abstract}
This paper proposes a strategy to solve the problems of the conventional s-version of finite element method (SFEM) fundamentally. 
Because SFEM can reasonably model an analytical domain by superimposing meshes with different spatial resolutions, it has intrinsic advantages of local high accuracy, low computation time, and simple meshing procedure.
However, it has disadvantages such as accuracy of numerical integration and matrix singularity. Although several additional techniques have been proposed to mitigate these limitations, they are computationally expensive or ad-hoc, and detract from the method’s strengths.  
To solve these issues, we propose a novel strategy called B-spline based SFEM.
To improve the accuracy of numerical integration, we employed cubic B-spline basis functions with $C^2$-continuity across element boundaries as the global basis functions. 
To avoid matrix singularity, we applied different basis functions to different meshes. Specifically, we employed the Lagrange basis functions as local basis functions.
The numerical results indicate that using the proposed method, numerical integration can be calculated with sufficient accuracy without any additional techniques used in conventional SFEM.
Furthermore, the proposed method avoids matrix singularity and is superior to conventional methods in terms of convergence for solving linear equations.
Therefore, the proposed method has the potential to reduce computation time while maintaining a comparable accuracy to conventional SFEM.
\end{abstract}



\begin{keyword}
s-version of finite element method \sep mesh superposition method \sep B-spline basis functions \sep localized mesh refinement


\end{keyword}

\end{frontmatter}


\section{Introduction}
The modeling of flow dynamics with moving boundaries and interfaces, including fluid--structure interactions \citep{ishihara2009two}, free-surface flows \citep{queutey2007interface}, and two-fluid flows \citep{qian2006free}, plays a prominent role in many scientific and engineering fields.
Depending on the nature of these problems, we can use an interface-tracking or interface-capturing method for their computation.

In an interface-tracking method, such as arbitrary Lagrangian--Eulerian (ALE) schemes \citep{hirt1974arbitrary} and deforming-spatial-domain/stabilized space--time (DSD/SST) \citep{tezduyar1991stabilized,tezduyar1992computation}, as the interfaces move and the fluid domain changes its shape, the mesh moves to adjust to the shape change and follow the interfaces. 
Moving the fluid mesh to follow the interfaces enables us to control the mesh resolution across the entire domain, produce a high-resolution representation of the boundary layers, and obtain highly accurate solutions in such critical flow regions. As we move the mesh, if the element distortion exceeds the threshold for good accuracy, a remeshing (i.e., mesh regenerating) is performed.
Although some advanced mesh update methods that aim to decrease the frequency of remeshing and sustain the good quality of elements near solid surfaces have been developed \citep{takizawa2020low,tonon2021linear}, these approaches incur additional non-negligible computation time and significant coding effort. In addition, even if an algorithm is robust to small interface movements, it may exhibit numerical instabilities when dealing with large deformations, movements, and contacts of interfaces \citep{sahin2009arbitrary}.

Generally, an interface-capturing method is used to solve these problems of interface-tracking methods. In this approach, flow fields are represented by a fixed Eulerian mesh regardless of the change in interface.
Although a fixed Eulerian mesh by itself cannot be used to model and simulate a complex moving geometry, combining it with interface representation approaches enables us to handle boundary conditions at interfaces.
Previously developed approaches include the immersed boundary (IB) methods \citep{peskin1972flow}, extended immersed boundary (EIB) method \citep{wang2004extended}, immersed finite element (IFE) method \citep{zhang2004immersed}, distributed Lagrange multiplier/fictitious domain (DLM/FD) method \citep{glowinski1999distributed}, extended FEM (XFEM) \citep{wagner2001extended}, finite cover method (FCM) \citep{terada2003finite}, cut-cell methods with marker particles \citep{udaykumar1996elafint}, and level-set methods \citep{dunne2006eulerian}.
Owing to their advantages associated with mesh processing and robustness against moving boundaries, interface-capturing approaches have been employed in various applications that involve complex geometries and large boundary movements ranging from incompressible flows to turbulence flows \citep{kan2021numerical} and FSI phenomena \citep{souza2022multi,kawakami2022fluid}. More details pertaining to these approaches can be found in several reviews \citep{osher2001level,kim2019immersed,huang2019recent}.
However, it remains difficult to achieve locally high resolution using interface-capturing approaches. 
To obtain highly accurate solutions in critical flow domains such as boundary layers, localized fine meshes are required. If a uniform mesh is used, this requirement is inevitably extended to the entire computational domain, and the resulting mesh may exceed storage capacity.
For this purpose, several adaptive mesh refinement schemes \citep{roma1999adaptive,hartmann2008adaptive,griffith2012immersed,salih2019thin,borker2019mesh,aldlemy2020adaptive}, overlapping schemes \citep{henshaw2008parallel,massing2014stabilized,bathe2017finite,huang2021convergence}, and hybrid schemes that merge concepts from capturing methods and ALE formulations \citep{gerstenberger2008enhancement} have been proposed.
However, these refinement algorithms often fail to guarantee cell conformity and consistent interpolation of the adapted meshes, or are highly complex.

By contrast, \citet{fish1992s} proposed the s-version of finite element method (SFEM), another localized mesh refinement approach. SFEM uses two-level FEM meshes - a global mesh and a local mesh - to model the target domain, where a fine local mesh(es) representing local features is superposed on the relatively coarse global mesh that represents the entire analytical domain.
Variables in the mesh superposing region are given by the sum of those in the global and local meshes.
Note that the local mesh can be inserted into an arbitrary part of the domain independently of the global mesh, so that the complex meshing procedure can be avoided.
SFEM has been successfully applied to various engineering problems, such as 
stress analyses of laminated composites \citep{fish1992multi,FISH1993363,REDDY199321,FISH1994135,ANGIONI2011780,ANGIONI2012559,Chen2014-ro,jiao2015adaptive,Jiao2015-fo,KUMAGAI2017136,Sakata2020-dq}, 
mesoscopic analyses of particulate composites\citep{Okada2004-el,Okada2004-co}, 
multiscale analyses of fibre-reinforced composites \citep{vorobiov2017mesh,sakata2022mesh}, 
concurrent multiscaling \citep{fish1993multiscale,sun2018variant,cheng2022multiscale}, 
multiscale analyses of porous materials \citep{takano2003multi,takano2004three,kawagai2006image,tsukino2015multiscale}, 
dynamic analyses of transient problems \citep{yue2005adaptive,yue2007adaptive}, and
shape and topology optimization problems \citep{wang2006moving}.
Fracture mechanics problems, such as fatigue crack and dynamic crack propagation, are also major applications of SFEM \citep{fish1993adaptive,FISH1994135,lee2004combined,okada2005fracture,okada2007application,fan2008rs,nakasumi2008crack,kikuchi2012crack,kikuchi2014fatigue,wada2014fatigue,kikuchi2016crack,xu2018study,kishi2020dynamic,cheng2022multiscale,he2023strategy,cheng2023application}. Furthermore, SFEM has been combined with XFEM \citep{lee2004combined,nakasumi2008crack,ANGIONI2011780,ANGIONI2012559,Jiao2015-fo} and phase field modeling \citep{cheng2022multiscale,cheng2023application}.
The approximation concept used in SFEM has been applied to the coupling of peridynamics and FEM \citep{sun2019superposition,sun2022numerical,sun2022parallel,sun2023stabilized}.

Despite its advantages, SFEM presents two challenges.
The first challenge is the inaccuracy of numerical integration based on Gaussian quadrature, occurring when the global and local elements exhibit partial mutual superposition.
The Lagrange basis functions used in the conventional SFEM have $C^0$-continuity, and their derivatives have discontinuities across element boundaries.
Owing to their low continuity, if the integral domain contains the Lagrange element boundaries, the integrands are often discontinuous. As a result, the accuracy of the Gaussian quadrature degrades.
To improve the accuracy, some approaches subdivide the integral domain into several subdomains \citep{FISH1993363,FISH1994135,Okada2004-el,Okada2004-co,okada2005fracture,okada2007application,he2023strategy}, while others apply high-order Gaussian quadrature \citep{fish1992s,lee2004combined,nakasumi2008crack,kishi2020dynamic}. 
However, all of them require large computation time.
The second challenge inherent to SFEM is the singularity of the matrix.
In the SFEM framework, two or more finite meshes are superposed, and the basis functions in said meshes are not guaranteed to be linearly independent of each other.
If the basis functions in one mesh can be represented as a linear combination of those in other meshes, matrix singularity occurs.
\citet{ooya2009linear} pointed out that once the problem arises, a linear equation solver using an iterative method, such as the conjugate gradient method, either fails to converge or is extremely slow to converge to the solution.
Although some approaches have been proposed to solve this problem \citep{fish1992s,FISH1994135,ANGIONI2011780,ANGIONI2012559,yue2005adaptive,yue2007adaptive,fan2008rs,nakasumi2008crack,park2003efficient,ooya2009linear}, many of them are ad-hoc methods that represent slight modifications to the models, or incur additional computation time.
More details on the difficulties of conventional SFEM are discussed in Section \ref{sec:difficulties_in_lagrange}.

In this study, we propose a new SFEM framework to avoid these problems fundamentally.
To improve the accuracy of numerical integration, functions with high continuity across the element boundaries were used as global basis functions such that the integrands are smooth and continuous.
Specifically, we applied cubic B-spline basis functions, which have $C^2$-continuity across the element boundaries, as the global basis functions.
Furthermore, to address matrix singularity, we applied different types of functions as the basis functions in different meshes.
Specifically, we employed Lagrange basis functions as the local basis functions. 
Note that unlike the conventional method, our framework does not require additional and computationally expensive techniques to address these issues.
Thus, the proposed method is expected to achieve the same level of accuracy with less computation time than the conventional method.

The remainder of this paper is organized as follows.
Basic formulations and the concept of SFEM are presented in Section \ref{sec:formulation_of_s-fem}.
The formulations and difficulties of conventional SFEM, and the formulations and the advantages of our proposed B-spline based SFEM method, are introduced in Section \ref{sec:bsfem}. 
The proposed method is verified in Section \ref{sec:verification}.
Finally, the conclusions of this study are presented in Section \ref{sec:conclusions}.

\section{Basic formulation of SFEM}
\label{sec:formulation_of_s-fem}
This section reviews the basic formulations and underlying concept of SFEM.

The target problem of this study is Poisson's equation with Dirichlet boundary conditions given as follows:
\begin{subequations}\label{eq:poisson}
    \begin{align}
        \Delta{u} + f &= 0 \quad \mathrm{in} \; \Omega, \\
        u &= g \quad \mathrm{on} \; \Gamma_D, \label{eq:dirichlet}
    \end{align}
\end{subequations}
where $\Omega$ is the domain and $\Gamma$ is its boundary, consisting of $\Gamma_D = \Gamma$. The function $u : \Omega \rightarrow \mathbb{R}$ is a trial solution and the function $f : \Omega \rightarrow \mathbb{R}$ is given. Eq. \eqref{eq:dirichlet} represents the Dirichlet boundary conditions.

We define the trial solution space $\mathcal{S}$ and test function space $\mathcal{V}$ as
\begin{equation}
    \mathcal{S} = \{ u \mid u \in H^1\left( \Omega \right), u|_{\Gamma_D} = g \}
\end{equation}
and
\begin{equation}
    \mathcal{V} = \{ w \mid w \in H^1\left( \Omega \right), w|_{\Gamma_D} = 0 \},
\end{equation}
respectively, where $H^1\left( \Omega \right)$ is the Sobolev space.

The resulting weak form of the problem is:
Given $f$ and $g$, find $u \in \mathcal{S}$ such that for all $w \in \mathcal{V}$,
\begin{equation}
    a_{\Omega} \left( w, u \right) = L_{\Omega} \left( w \right),
\end{equation}
where
\begin{equation}
    a_{\Omega} \left( w, u \right) = \int_{\Omega} \nabla{w} \cdot \nabla{u} \; d \Omega,
\end{equation}
and
\begin{equation}
    L_{\Omega} \left( w \right) = \int_{\Omega} w f \; d \Omega.
\end{equation}
Here, $a_{\Omega} \left( \cdot, \cdot \right)$ is a bilinear form, and $L_{\Omega} \left( \cdot \right)$ is a linear functional.

\begin{figure}[t]
	\centering
	\includegraphics[bb=0 0 993 564,width=9cm]{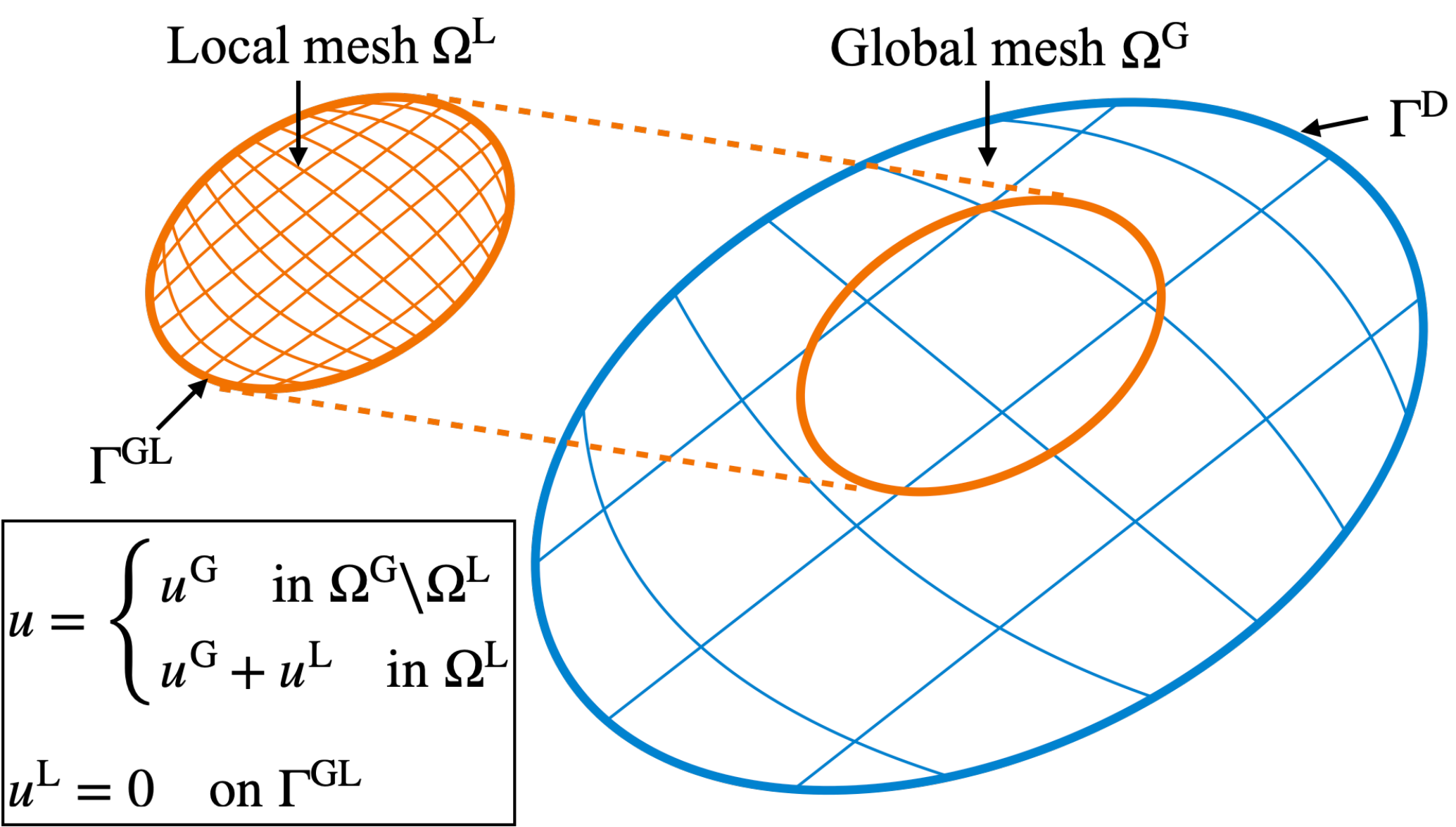}
	\caption{Global and local meshes defined in SFEM.}
	\label{fig:sfem_meshes}
\end{figure}

In the framework of SFEM \citep{fish1992s}, the domain $\Omega$ is discretized by some finite element meshes defined with mutual independence. In many cases, as shown in Figure \ref{fig:sfem_meshes}, one relatively coarse mesh may represent the entire domain $\Omega^{\mathrm{G}}$ that corresponds to domain $\Omega$, known as the global mesh. The local domain $\Omega^{\mathrm{L}}$ requiring high resolution is discretized by the finer local mesh, which is superimposed on the global mesh. 
The local domain $\Omega^{\mathrm{L}}$ is assumed to be included in the global domain $\Omega^{\mathrm{G}}$ as $\Omega^{\mathrm{L}}$ $\subseteq$ $\Omega^{\mathrm{G}}$.
In the present work, although the local domain $\Omega^{\mathrm{L}}$ must not extend outside of the global domain $\Omega^{\mathrm{G}}$, the boundary of the local domain $\Omega^{\mathrm{L}}$ is permitted to overlap with the boundary of the global domain $\Omega^{\mathrm{G}}$.
Because meshes are defined independently from each other without considering their mutual consistency, mesh generation can be simplified.
Subscripts $\mathrm{G}$ and $\mathrm{L}$ respectively represent the quantities of the global and local meshes, as illustrated in Figure \ref{fig:sfem_meshes}.

In the formulation of SFEM, the trial solution $u$ in different regions is defined as
\begin{align} 
    u = 
    \begin{cases}
        {u^{\mathrm{G}} \quad \mathrm{in} \; \Omega^{\mathrm{G}} \setminus \Omega^{\mathrm{L}}},\\
        {u^{\mathrm{G}} + u^{\mathrm{L}} \quad \mathrm{in} \; \Omega^{\mathrm{L}}}.
    \end{cases}
\end{align}
To ensure continuity between the global and local meshes, the following Dirichlet boundary condition is imposed:
\begin{equation}
    u^{\mathrm{L}} = 0 \quad \mathrm{on} \; \Gamma^{\mathrm{GL}},
    \label{eq:dirichlet_local}
\end{equation}
where $\Gamma^{\mathrm{GL}}$ is the boundary of the local domain.
Because the global domain $\Omega^{\mathrm{G}}$ corresponds to domain $\Omega$, the following Dirichlet boundary condition is imposed:
\begin{equation}
    u^{\mathrm{G}} = g \quad \mathrm{on} \; \Gamma_D.
\end{equation}
Based on the Galerkin method, the test function is defined as
\begin{align} 
    w = 
    \begin{cases}
        {w^{\mathrm{G}} \quad \mathrm{in} \; \Omega^{\mathrm{G}} \setminus \Omega^{\mathrm{L}}},\\
        {w^{\mathrm{G}} + w^{\mathrm{L}} \quad \mathrm{in} \; \Omega^{\mathrm{L}}}.
    \end{cases}
\end{align}
We therefore define the trial solution spaces $\mathcal{S}^{\mathrm{G}}$, $\mathcal{S}^{\mathrm{L}}$ and test function spaces $\mathcal{V}^{\mathrm{G}}$, $\mathcal{V}^{\mathrm{L}}$ as
\begin{align}
    \mathcal{S}^{\mathrm{G}} &= \{ u^{\mathrm{G}} \mid u^{\mathrm{G}} \in H^1\left( \Omega^{\mathrm{G}} \right), u^{\mathrm{G}}|_{\Gamma_D} = g \},\\
    \mathcal{S}^{\mathrm{L}} &= \{ u^{\mathrm{L}} \mid u^{\mathrm{L}} \in H^1\left( \Omega^{\mathrm{L}} \right), u^{\mathrm{L}}|_{\Gamma^{\mathrm{GL}}} = 0 \},\\
    \mathcal{V}^{\mathrm{G}} &= \{ w^{\mathrm{G}} \mid w^{\mathrm{G}} \in H^1\left( \Omega^{\mathrm{G}} \right), w^{\mathrm{G}}|_{\Gamma_D} = 0 \},\\
    \mathcal{V}^{\mathrm{L}} &= \{ w^{\mathrm{L}} \mid w^{\mathrm{L}} \in H^1\left( \Omega^{\mathrm{L}} \right), w^{\mathrm{L}}|_{\Gamma^{\mathrm{GL}}} = 0 \}.
\end{align}

The resulting weak form of Eq. \eqref{eq:poisson} in the SFEM formulation is: Given $f$ and $g$, find $\left( u^{\mathrm{G}}, u^{\mathrm{L}} \right) \in \mathcal{S}^{\mathrm{G}} \times \mathcal{S}^{\mathrm{L}}$  such that for all $\left( w^{\mathrm{G}}, w^{\mathrm{L}} \right) \in \mathcal{V}^{\mathrm{G}} \times \mathcal{V}^{\mathrm{L}}$
\begin{equation}
    a_{\Omega}^{'} \left( w, u \right) = L_{\Omega}^{'} \left( w \right),
    \label{eq:variational_form_sfem_1}
\end{equation}
where
\begin{align} 
    a_{\Omega}^{'} \left( w, u \right) =& \:
    a_{\Omega^{\mathrm{G}} \setminus \Omega^{\mathrm{L}}} \left( w^{\mathrm{G}}, u^{\mathrm{G}} \right) + a_{\Omega^{\mathrm{L}}} \left( w^{\mathrm{G}} + w^{\mathrm{L}}, u^{\mathrm{G}} + u^{\mathrm{L}} \right),\label{eq:sfem_a_1}\\
    L_{\Omega}^{'} \left( w \right) =& \:
    L_{\Omega^{\mathrm{G}} \setminus \Omega^{\mathrm{L}}} \left( w^{\mathrm{G}} \right) + L_{\Omega^{\mathrm{L}}} \left( w^{\mathrm{G}} + w^{\mathrm{L}} \right).\label{eq:sfem_L_1}
\end{align}
Recalling the bilinearity of $a_{\Omega} \left( \cdot, \cdot \right)$ and linearity of $L_{\Omega} \left( \cdot \right)$, we can rewrite Eqs. \eqref{eq:sfem_a_1} and \eqref{eq:sfem_L_1} as
\begin{equation}
  \begin{split} 
    a_{\Omega}^{'} \left( w, u \right)
    = \: &a_{\Omega^{\mathrm{G}}} \left( w^{\mathrm{G}}, u^{\mathrm{G}} \right)
    + a_{\Omega^{\mathrm{L}}} \left( w^{\mathrm{G}}, u^{\mathrm{L}} \right) \\
    &+ a_{\Omega^{\mathrm{L}}} \left( w^{\mathrm{L}}, u^{\mathrm{G}} \right)
    + a_{\Omega^{\mathrm{L}}} \left( w^{\mathrm{L}}, u^{\mathrm{L}} \right),
  \end{split}
\end{equation}
\begin{align}
    L' \left( w \right) =
    L_{\Omega^{\mathrm{G}}} \left( w^{\mathrm{G}} \right) +
    L_{\Omega^{\mathrm{L}}} \left( w^{\mathrm{L}} \right).
\end{align}

To convert this weak problem statement into a coupled system of linear algebraic equations, we apply Galerkin's method and work in finite-dimensional subspaces $\left( \mathcal{S}^{\mathrm{G}} \right)^h \subset \mathcal{S}^{\mathrm{G}}$, $\left( \mathcal{S}^{\mathrm{L}} \right)^h \subset \mathcal{S}^{\mathrm{L}}$, $\left( \mathcal{V}^{\mathrm{G}} \right)^h \subset \mathcal{V}^{\mathrm{G}}$, and $\left( \mathcal{V}^{\mathrm{L}} \right)^h \subset \mathcal{V}^{\mathrm{L}}$. 

We have a given function $g^h \in \left( \mathcal{S}^{\mathrm{G}} \right)^h$ such that $g^h|_{\Gamma_D} = g$, and thus for every $\left( u^{\mathrm{G}} \right)^h \in \left( \mathcal{S}^{\mathrm{G}} \right)^h$ we have a unique decomposition
\begin{equation}
    \left( u^{\mathrm{G}} \right)^h = \left( v^{\mathrm{G}} \right)^h + g^h,
\end{equation}
where $\left( v^{\mathrm{G}} \right)^h \in \left( \mathcal{V}^{\mathrm{G}} \right)^h$. 

Therefore the Galerkin approximation of Eq. \eqref{eq:variational_form_sfem_1} is:
Find $\left( u^{\mathrm{G}} \right)^h = \left( v^{\mathrm{G}} \right)^h + g^h$ and $\left( u^{\mathrm{L}} \right)^h = \left( v^{\mathrm{L}} \right)^h$, where $\left( \left( v^{\mathrm{G}} \right)^h, \left( v^{\mathrm{L}} \right)^h \right) \in \left( \mathcal{V}^{\mathrm{G}} \right)^h \times \left( \mathcal{V}^{\mathrm{L}} \right)^h$, such that for all $\left( \left( w^{\mathrm{G}} \right)^h, \left( w^{\mathrm{L}} \right)^h \right) \in \left( \mathcal{V}^{\mathrm{G}} \right)^h \times \left( \mathcal{V}^{\mathrm{L}} \right)^h$
\begin{equation}
    a_{\Omega}^{'} \left( w^h, v^h \right) = L_{\Omega}^{'} \left( w^h \right) - a_{\Omega}^{'} \left( w^h, g^h \right),
    \label{eq:variational_form_sfem_2}
\end{equation}
where
\begin{equation}
  \begin{split} 
    a_{\Omega}^{'} \left( w^h, v^h \right)
    = \: 
    &a_{\Omega^{\mathrm{G}}} \left( \left( w^{\mathrm{G}} \right)^h, \left( v^{\mathrm{G}} \right)^h \right)
    + a_{\Omega^{\mathrm{L}}} \left( \left( w^{\mathrm{G}} \right)^h, \left( v^{\mathrm{L}} \right)^h \right) \\
    &+ a_{\Omega^{\mathrm{L}}} \left( \left( w^{\mathrm{L}} \right)^h, \left( v^{\mathrm{G}} \right)^h \right)
    + a_{\Omega^{\mathrm{L}}} \left( \left( w^{\mathrm{L}} \right)^h, \left( v^{\mathrm{L}} \right)^h \right),
  \end{split}
\end{equation}
\begin{align}
    L' \left( w \right) =
    L_{\Omega^{\mathrm{G}}} \left( \left( w^{\mathrm{G}} \right)^h \right) +
    L_{\Omega^{\mathrm{L}}} \left( \left( w^{\mathrm{L}} \right)^h \right),
\end{align}
and
\begin{align}
    a_{\Omega}^{'} \left( w^h, g^h \right) =
    a_{\Omega^{\mathrm{G}}} \left( \left( w^{\mathrm{G}} \right)^h, g^h \right) +
    a_{\Omega^{\mathrm{G}}} \left( \left( w^{\mathrm{L}} \right)^h, g^h \right).
\end{align}

We define $\boldsymbol{\eta}^{\mathrm{G}}$ and $\boldsymbol{\eta}^{\mathrm{L}}$ to be the sets containing the indices of all global basis functions $N_A^{\mathrm{G}}, A \in \boldsymbol{\eta}^{\mathrm{G}}$ and local basis functions $N_A^{\mathrm{L}}, A \in \boldsymbol{\eta}^{\mathrm{L}}$, respectively.
Similarly, we let $\boldsymbol{\eta}_D^{\mathrm{G}} \subset \boldsymbol{\eta}^{\mathrm{G}}$ be the set containing the indices of all of global basis functions that are non-zero on $\Gamma_D$. Thus, $\left( u^{\mathrm{G}} \right)^h$ and $\left( u^{\mathrm{L}} \right)^h$ can be expressed as
\begin{align}
    \left( u^{\mathrm{G}} \right)^h &= 
    \sum_{A \in \boldsymbol{\eta}^{\mathrm{G}} - \boldsymbol{\eta}_D^{\mathrm{G}}} N_A^{\mathrm{G}} d_A^{\mathrm{G}} + 
    \sum_{B \in \boldsymbol{\eta}_D^{\mathrm{G}}} N_B^{\mathrm{G}} g_B = 
    \sum_{A \in \boldsymbol{\eta}^{\mathrm{G}} - \boldsymbol{\eta}_D^{\mathrm{G}}} N_A^{\mathrm{G}} d_A^{\mathrm{G}} + g^h,
    \label{eq:interporation_uG} \\
    \left( u^{\mathrm{L}} \right)^h &= 
    \sum_{A \in \boldsymbol{\eta}^{\mathrm{L}}} N_A^{\mathrm{L}} d_A^{\mathrm{L}}. 
    \label{eq:interporation_uL}
\end{align}

Similarly, $\left( w^{\mathrm{G}} \right)^h$ and $\left( w^{\mathrm{L}} \right)^h$ can be expressed as
\begin{align}
    \left( w^{\mathrm{G}} \right)^h &= \sum_{A \in \boldsymbol{\eta}^{\mathrm{G}} - \boldsymbol{\eta}_D^{\mathrm{G}}} N_A^{\mathrm{G}} e_A^{\mathrm{G}}, \label{eq:interporation_wG} \\
    \left( w^{\mathrm{L}} \right)^h &= \sum_{A \in \boldsymbol{\eta}^{\mathrm{L}}} N_A^{\mathrm{L}} e_A^{\mathrm{L}}. \label{eq:interporation_wL}
\end{align}

Inserting Eqs. \eqref{eq:interporation_uG}, \eqref{eq:interporation_uL}, \eqref{eq:interporation_wG}, and \eqref{eq:interporation_wL} into Eq. \eqref{eq:variational_form_sfem_2} yields
\begin{equation}
  \begin{split} 
    \sum_{A \in \boldsymbol{\eta}^{\mathrm{G}} - \boldsymbol{\eta}_D^{\mathrm{G}}} e_A^{\mathrm{G}} \Biggl(
    \sum_{B \in \boldsymbol{\eta}^{\mathrm{G}} - \boldsymbol{\eta}_D^{\mathrm{G}}} 
    a_{\Omega^{\mathrm{G}}} \left( N_A^{\mathrm{G}}, N_B^{\mathrm{G}} \right)d_B^{\mathrm{G}} +
    \sum_{C \in \boldsymbol{\eta}^{\mathrm{L}}} 
    a_{\Omega^{\mathrm{L}}} \left( N_A^{\mathrm{G}}, N_C^{\mathrm{L}} \right)d_C^{\mathrm{L}} \Biggr. \\
    \Biggl. - L_{\Omega^{\mathrm{G}}} \left( N_A^{\mathrm{G}} \right)
    + a_{\Omega^{\mathrm{G}}} \left( N_A^{\mathrm{G}}, g^h \right)
    \Biggr) \\
    + \sum_{D \in \boldsymbol{\eta}^{\mathrm{L}}} e_D^{\mathrm{L}} \Biggl(
    \sum_{B \in \boldsymbol{\eta}^{\mathrm{G}} - \boldsymbol{\eta}_D^{\mathrm{G}}} 
    a_{\Omega^{\mathrm{L}}} \left( N_D^{\mathrm{L}}, N_B^{\mathrm{G}} \right)d_B^{\mathrm{G}} +
    \sum_{C \in \boldsymbol{\eta}^{\mathrm{L}}} 
    a_{\Omega^{\mathrm{L}}} \left( N_D^{\mathrm{L}}, N_C^{\mathrm{L}} \right)d_C^{\mathrm{L}} \Biggr. \\
    \Biggl. - L_{\Omega^{\mathrm{L}}} \left( N_D^{\mathrm{L}} \right)
    + a_{\Omega^{\mathrm{G}}} \left( N_D^{\mathrm{L}}, g^h \right)
    \Biggr) = 0.
  \end{split}
\end{equation}
As the $e_A^{\mathrm{G}}$'s and $e_D^{\mathrm{L}}$'s are arbitrary, the terms in parentheses must be identically zero. Thus, for $A \in \boldsymbol{\eta}^{\mathrm{G}} - \boldsymbol{\eta}_D^{\mathrm{G}}$ and $D \in \boldsymbol{\eta}^{\mathrm{L}}$,
\begin{multline}
    \sum_{B \in \boldsymbol{\eta}^{\mathrm{G}} - \boldsymbol{\eta}_D^{\mathrm{G}}} 
    a_{\Omega^{\mathrm{G}}} \left( N_A^{\mathrm{G}}, N_B^{\mathrm{G}} \right)d_B^{\mathrm{G}} +
    \sum_{C \in \boldsymbol{\eta}^{\mathrm{L}}} 
    a_{\Omega^{\mathrm{L}}} \left( N_A^{\mathrm{G}}, N_C^{\mathrm{L}} \right)d_C^{\mathrm{L}}\\
    = L_{\Omega^{\mathrm{G}}} \left( N_A^{\mathrm{G}} \right)
    - a_{\Omega^{\mathrm{G}}} \left( N_A^{\mathrm{G}}, g^h \right),
    \label{eq:befor_matrix1}
\end{multline}
\begin{multline}
    \sum_{B \in \boldsymbol{\eta}^{\mathrm{G}} - \boldsymbol{\eta}_D^{\mathrm{G}}} 
    a_{\Omega^{\mathrm{L}}} \left( N_D^{\mathrm{L}}, N_B^{\mathrm{G}} \right)d_B^{\mathrm{G}} +
    \sum_{C \in \boldsymbol{\eta}^{\mathrm{L}}} 
    a_{\Omega^{\mathrm{L}}} \left( N_D^{\mathrm{L}}, N_C^{\mathrm{L}} \right)d_C^{\mathrm{L}}\\
    = L_{\Omega^{\mathrm{L}}} \left( N_D^{\mathrm{L}} \right)
    - a_{\Omega^{\mathrm{G}}} \left( N_D^{\mathrm{L}}, g^h \right).
    \label{eq:befor_matrix2}
\end{multline}
Proceeding to define
\begin{align}
    K_{AB} &= a_{\Omega^{\mathrm{G}}} \left( N_A^{\mathrm{G}}, N_B^{\mathrm{G}} \right), \\
    K_{AC} &= a_{\Omega^{\mathrm{L}}} \left( N_A^{\mathrm{G}}, N_C^{\mathrm{L}} \right), \\
    K_{DB} &= a_{\Omega^{\mathrm{L}}} \left( N_D^{\mathrm{L}}, N_B^{\mathrm{G}} \right), \\
    K_{DC} &= a_{\Omega^{\mathrm{L}}} \left( N_D^{\mathrm{L}}, N_C^{\mathrm{L}} \right), \\
    F_A &= L_{\Omega^{\mathrm{G}}} \left( N_A^{\mathrm{G}} \right)
    - a_{\Omega^{\mathrm{G}}} \left( N_A^{\mathrm{G}}, g^h \right),\\
    F_D &= L_{\Omega^{\mathrm{L}}} \left( N_D^{\mathrm{L}} \right)
    - a_{\Omega^{\mathrm{G}}} \left( N_D^{\mathrm{L}}, g^h \right),\\
    \boldsymbol{K}^{\mathrm{GG}} &=  \lbrack K_{AB} \rbrack, \\
    \boldsymbol{K}^{\mathrm{GL}} &=  \lbrack K_{AC} \rbrack, \\
    \boldsymbol{K}^{\mathrm{LG}} &=  \lbrack K_{DB} \rbrack, \\
    \boldsymbol{K}^{\mathrm{LL}} &=  \lbrack K_{DC} \rbrack, \\
    \boldsymbol{F}^{\mathrm{G}} &= \lbrace F_A \rbrace, \\
    \boldsymbol{F}^{\mathrm{L}} &= \lbrace F_D \rbrace, \\
    \boldsymbol{d}^{\mathrm{G}} &= \lbrace d_B^{\mathrm{G}} \rbrace, \\
    \boldsymbol{d}^{\mathrm{L}} &= \lbrace d_C^{\mathrm{L}} \rbrace,
\end{align}
and
\begin{align}
    \boldsymbol{K} &= 
    \begin{bmatrix}
        \boldsymbol{K}^{\mathrm{GG}} & \boldsymbol{K}^{\mathrm{GL}} \\
        \boldsymbol{K}^{\mathrm{LG}} & \boldsymbol{K}^{\mathrm{LL}} \\
    \end{bmatrix}, \label{eq:submatrix_K} \\
    \boldsymbol{F} &=
    \begin{bmatrix}
        \boldsymbol{F}^{\mathrm{G}} \\
        \boldsymbol{F}^{\mathrm{L}} \\
    \end{bmatrix}, \\
    \boldsymbol{d} &=
    \begin{bmatrix}
        \boldsymbol{d}^{\mathrm{G}} \\
        \boldsymbol{d}^{\mathrm{L}} \\
    \end{bmatrix},
\end{align}
for $A, B \in \boldsymbol{\eta}^{\mathrm{G}} - \boldsymbol{\eta}_D^{\mathrm{G}}$ and $C, D \in \boldsymbol{\eta}^{\mathrm{L}}$, we can rewrite Eqs. \eqref{eq:befor_matrix1} and \eqref{eq:befor_matrix2} as simultaneous linear equations:
\begin{align}
    \boldsymbol{K} \boldsymbol{d} = \boldsymbol{F},
    \label{eq:system_sfem}
\end{align}
where $\boldsymbol{K}^{\mathrm{GG}}$ and $\boldsymbol{K}^{\mathrm{LL}}$ are the submatrices for the global and local meshes, respectively, and $\boldsymbol{K}^{\mathrm{GL}}$ and $\boldsymbol{K}^{\mathrm{LG}}$ are submatrices representing the relationship between said meshes.

The primary challenges of standard SFEM are the difficulty in exact integration of the submatrices $\boldsymbol{K}^{\mathrm{GL}}$ and $\boldsymbol{K}^{\mathrm{LG}}$, and the singularity of the matrix $\boldsymbol{K}$. These topics are discussed in Section \ref{sec:difficulties_in_lagrange}.

\section{B-spline based s-version of finite element method (BSFEM)}
\label{sec:bsfem}
This section introduces the concept and formulations of our proposed method.
Lagrange basis functions, and conventional SFEM problems using said functions, are defined in Sections \ref{sec:formulation_of_lagrange} and \ref{sec:difficulties_in_lagrange}, respectively.
Subsequently, B-spline basis functions and details of the proposed method are presented in Sections \ref{sec:formulation_of_b-spline} and \ref{sec:details_bsfem}, respectively.

\subsection{Formulation of Lagrange basis functions}
\label{sec:formulation_of_lagrange}
In this section, we briefly summarize the basics of Lagrange basis functions used in conventional SFEM.

Let us denote coordinates in the physical space and parent element by $\boldsymbol{x}$ and $\hat{\boldsymbol{\xi}}$, respectively.
To use Gaussian quadrature for integration, the interval of the parent element is defined as $[-1, 1]^d$, where $d$ is the space dimension.

In the parent element,
the Lagrange interpolation formula in one dimension is given by
\begin{equation}
    l_i^p \left( \hat{\xi} \right) = \prod_{j=1, j \neq i}^{p+1} \frac{\hat{\xi} - \hat{\xi}_j}{\hat{\xi}_i - \hat{\xi}_j},
\end{equation}
where $i=1, 2, \dotsc, p+1$, $p$ is the order of the polynomial, and the number of functions (or nodes) used in the parent element is $n = p+1$. 
Let $\hat{\xi}_i$ be the parent coordinate of node $i$ and $\hat{\xi}_1 = -1, \hat{\xi}_{p+1} = 1$.

The $p$th order Lagrange basis functions consist of the $p$th order Lagrange polynomials $l_i^p$.

\begin{figure}[t]
	\centering
	\includegraphics[bb=0 0 1057 388,width=12cm]{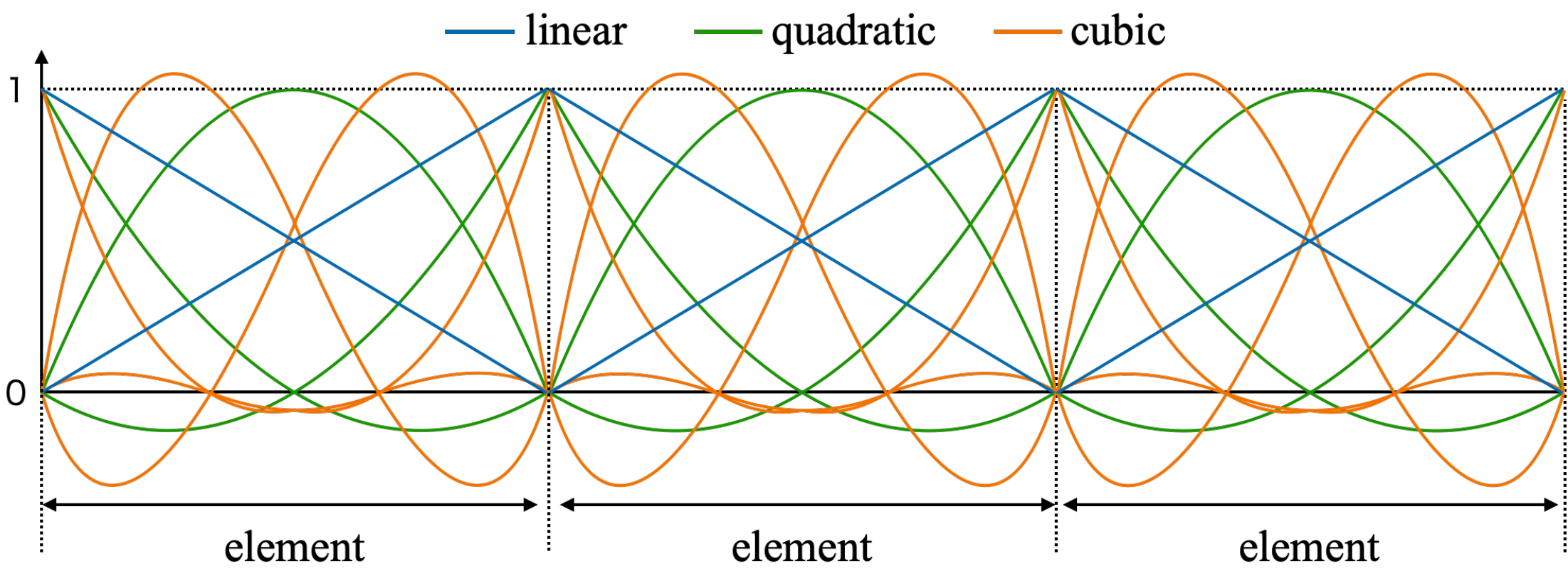}
	\caption{Linear, quadratic and cubic Lagrange basis functions.}
	\label{fig:lagrange_function}
\end{figure}

As shown in Figure \ref{fig:lagrange_function}, Lagrange basis functions have several essential features.
The most important feature noted in this study is that although each Lagrange basis function of order $p$ has $C^{p-1}$-continuous derivatives inside each element, it has $C^0$-continuity and its derivatives have discontinuity across their respective element boundaries.

Furthermore, each Lagrange basis function satisfies the interpolation property; that is,
\begin{equation}
    N_i \left( \hat{\xi}_j \right) = \delta_{ij}.
\end{equation}

\subsection{Difficulties in conventional SFEM using Lagrange basis functions}
\label{sec:difficulties_in_lagrange}
Two challenges with the conventional Lagrange-based SFEM are described in this section.
The first is the inaccuracy of numerical integration based on Gaussian quadrature.
The second is the singularity of the matrix, occurring if the basis functions in one mesh can be represented as a linear combination of the basis functions in other meshes.

\subsubsection{Inaccurate numerical integration for discontinuous functions}
\begin{figure}[t]
	\centering
	\includegraphics[bb=0 0 693 428,width=8cm]{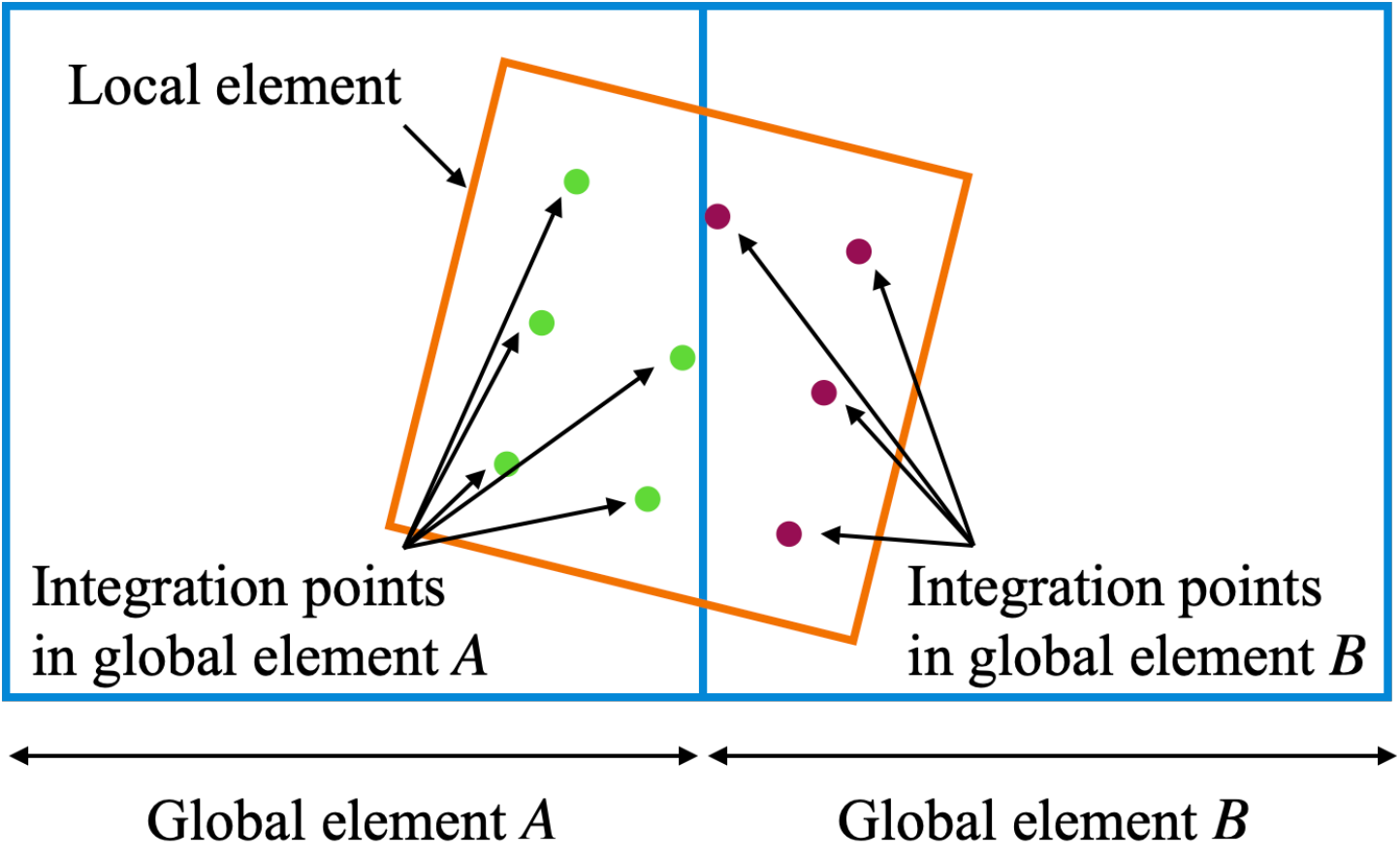}
	\caption{An example of mesh superimposition that causes inaccurate quadrature.}
	\label{fig:discontinuous_integration}
\end{figure}

The first challenge of conventional SFEM involves the numerical integration of the submatrices $\boldsymbol{K}^{\mathrm{GL}}$ and $\boldsymbol{K}^{\mathrm{LG}}$ in Eq. \eqref{eq:submatrix_K}.
The accuracy of numerical integration deteriorates based on Gaussian quadrature when a local element is superimposed on several global elements as shown in Figure \ref{fig:discontinuous_integration}.
In many cases, the integration of $\boldsymbol{K}^{\mathrm{GL}}$ and $\boldsymbol{K}^{\mathrm{LG}}$ is conducted by each element on the local mesh. 
Figure \ref{fig:discontinuous_integration} shows that a local element contains two different global elements $A, B$ and the boundaries between them.
The integrands often contain the first derivatives of global basis functions, which are Lagrange basis functions in conventional SFEM.
As mentioned in Section \ref{sec:formulation_of_lagrange}, the first derivatives of $p$th order Lagrange basis functions exhibit discontinuities across element boundaries. 
The resulting integrands often become discontinuous for local elements located at the boundaries of global elements.
Because the Gaussian quadrature scheme assumes that the integrands are smooth and continuous, the accuracy degrades.

\begin{figure}[t]
	\centering
	\includegraphics[bb=0 0 1185 292,width=13cm]{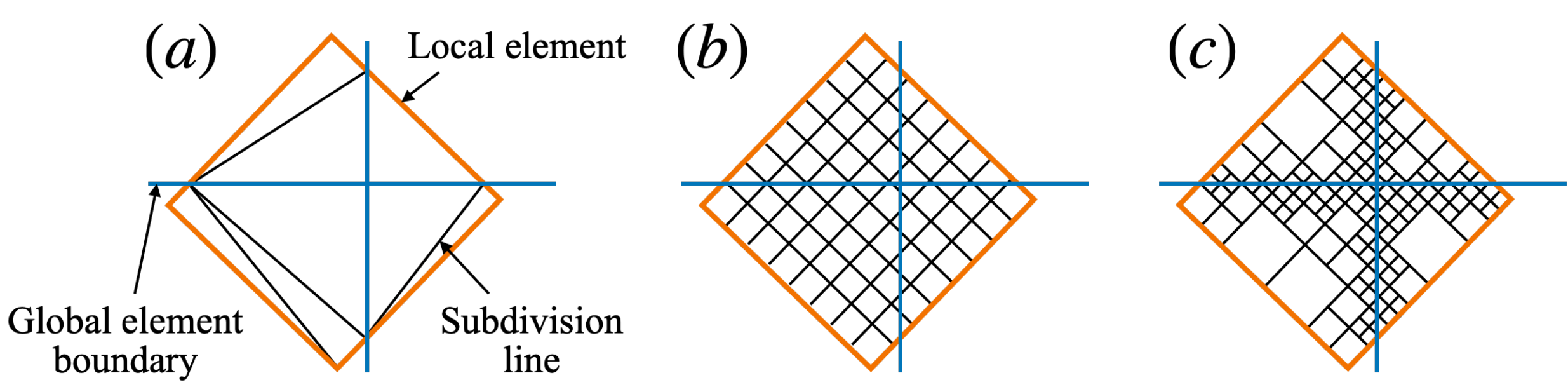}
	\caption{Examples of existing mesh-subdivision approaches for discontinuous integrands.}
	\label{fig:subdivision_technique}
\end{figure}

For the exact integration of such submatrices, \citet{FISH1993363} and \citet{FISH1994135} 
subdivided the local element into several subdomains, each of which corresponds to a global element, and performed Gaussian quadrature separately in each subdomain. (see Figure \ref{fig:subdivision_technique}(a))
However, this subdivision process requires a complex treatment of geometries and incurs substantial computation time.
Furthermore, the boundaries of each subdomain must be precisely defined, and the resulting subdomains may be arbitrary polygons that require further subdivision into simpler shapes.
The resulting complexity may be excessive when unstructured meshes are used.

To reduce the computation time other studies have applied the following approximate quadrature schemes.
As shown in Figure \ref{fig:subdivision_technique}(b), \citet{okada2007application} divided local elements into equal-sized square domains in the parent element coordinate space to confine numerical error due to discontinuous variations of the integrands in subdomains that contain discontinuities.
They also discussed the effect of this subdivision technique on accuracy \citep{okada2005fracture}.
In other studies \citep{Okada2004-el,Okada2004-co,he2023strategy}, local elements that contain the edges of the global mesh are divided recursively as shown in Figure \ref{fig:subdivision_technique}(c).
On the other hand, to integrate discontinuous functions without subdivision techniques, a high-order Gauss quadrature may be used \citep{fish1992s,lee2004combined,nakasumi2008crack,kishi2020dynamic,sawada2010high}.

Although more effective subdivision techniques have been devised, they still incur additional computation time, and the more complex the target problem, the greater the impact of deterioration in computational accuracy.

\subsubsection{Loss of solution uniqueness based on independency of basis functions}
\begin{figure}[t]
	\centering
	\includegraphics[bb=0 0 353 222,width=6cm]{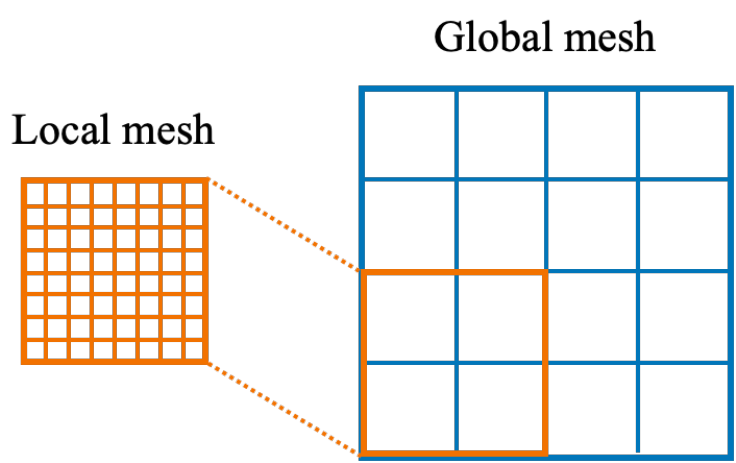}
	\caption{An example of mesh superimposition causing a loss of solution uniqueness.}
	\label{fig:loss_uniqueness}
\end{figure}

It is well-established that because the basis functions in two or more finite element meshes are not guaranteed to be linearly independent of each other, uniqueness of the decomposing numerical solution $u=u^{\mathrm{G}} + u^{\mathrm{L}}$ is sometimes lost \citep{fish1992s,FISH1994135,ANGIONI2011780,ANGIONI2012559,yue2005adaptive,yue2007adaptive,fan2008rs,nakasumi2008crack,park2003efficient,ooya2009linear}.
In other words, if the basis functions in one mesh can be represented as a linear combination of those in other meshes, a singularity of the matrix $\boldsymbol{K}$ in Eq. \eqref{eq:system_sfem} occurs. This is often the case if there is a patch of local elements having entire boundaries aligned along the global element sides as shown in Figure \ref{fig:loss_uniqueness}.
Once the problem arises, a linear equation solver using an iterative method, such as the conjugate gradient method, either fails to converge to the solution, or exhibits extremely slow convergence.

To solve this problem, 
several prior studies have employed structured meshes for both global and local meshes, suppressing the degrees of freedom at the nodes in local meshes that coincide with those at the global mesh to eliminate redundant degrees of freedom in the former \citep{park2003efficient,yue2005adaptive,yue2007adaptive,fan2008rs,ANGIONI2011780,ANGIONI2012559}.
Other approaches eliminate equations with zero (or close to zero) pivots encountered in the course of factorizing the equations corresponding to the local meshes, thereby ensuring rank sufficiency in the case of unstructured mesh superimposes \citep{FISH1994135}.
\citet{ooya2009linear} proposed an approach that systematically finds the linear dependencies of degrees of freedom, and suggested the possibility that the dependency results in an ill-conditioned matrix.
\citet{nakasumi2008crack} discretized the global mesh aslant.

However, many of these methods are ad-hoc approaches that slightly modify the underlying models, typically requiring additional computation time.
Consequently, approaches that constrain the degrees of freedom detract from low computation time and simplicity in the meshing procedure of SFEM.

\subsection{Formulation of B-spline basis functions}
\label{sec:formulation_of_b-spline}
The following section summarizes B-spline basis functions.
We note that these functions are employed to ensure a smooth global mesh discretization, which has significant benefits compared to Lagrange basis functions in the numerical integration of the submatrices $\boldsymbol{K}^{\mathrm{GL}}$ and $\boldsymbol{K}^{\mathrm{LG}}$.

We first describe the basic framework.
A knot vector in one dimension is a non-decreasing set of coordinates in the parameter space $\hat{\Omega}$ written as $\boldsymbol{\Xi} = \{ \xi_1, \xi_2, \dotsc, \xi_{n+p+1} \}$, where $\xi_i \in \mathbb{R}$ is the $i$th knot, $i$ is the knot index, $i=1, 2, \dotsc, n+p+1$, $p$ is the polynomial order, and $n$ is the number of B-spline basis functions. The knots partition the parametric space into elements.

For a given knot vector, the B-spline basis functions are defined recursively starting with piecewise constants $\left( p=0 \right)$
\begin{align} 
    N_{i,0} \left( \xi \right) = 
    \begin{cases}
        {1 \quad \mathrm{if} \; \xi_{i} \leq \xi < \xi_{i+1} },\\
        {0 \quad \mathrm{otherwise}}.
    \end{cases}
\end{align}

For $p=1, 2, 3, \dotsc, $ these functions are defined by
\begin{align}
    N_{i,p} \left( \xi \right) = 
    \frac{\xi - \xi_i}{\xi_{i+p} - \xi_i} N_{i,p-1}\left( \xi \right) +
    \frac{\xi_{i+p+1} - \xi}{\xi_{i+p+1} - \xi_{i+1}} N_{i+1,p-1}\left( \xi \right),
\end{align}
which is the Cox--de Boor recursion formula.

The derivatives of the functions are represented in terms of B-spline lower-order bases. For a given polynomial order $p$ and knot vector $\boldsymbol{\Xi}$, the derivative of the $i$th B-spline basis function is given by
\begin{align}
    \frac{d}{d \xi} N_{i,p} \left( \xi \right) = 
    \frac{p}{\xi_{i+p} - \xi_i} N_{i,p-1}\left( \xi \right) -
    \frac{p}{\xi_{i+p+1} - \xi_{i+1}} N_{i+1,p-1}\left( \xi \right).
\end{align}

B-spline curves are constructed by taking a linear combination of B-spline basis functions:
\begin{equation}
    \boldsymbol{C}\left( \xi \right) = \sum_{i=1}^{n} N_{i, p} \left( \xi \right) \boldsymbol{B}_i,
\end{equation}
where $\boldsymbol{B}_i \in \mathbb{R}^d$ are the control points.

In higher dimensions, B-spline basis functions are constructed from tensor products similarly to Lagrange basis functions. Given additional knot vectors $\boldsymbol{\mathcal{H}} = \{ \eta_1, \eta_2, \dotsc, \eta_{m+q+1} \}$ and $\boldsymbol{\mathcal{Z}} = \{ \zeta_1, \zeta_2, \dotsc, \zeta_{l+r+1} \}$, 
B-spline volumes are defined as
\begin{equation}
    \boldsymbol{V}\left( \xi, \eta, \zeta \right) = \sum_{i=1}^{n} \sum_{j=1}^{m} \sum_{k=1}^{l} N_{i, p}\left( \xi \right) M_{j, q} \left( \eta \right) L_{k, r} \left( \zeta \right) \boldsymbol{B}_{i, j, k},
\end{equation}
where $N_{i, p}\left( \xi \right)$, $M_{j, q} \left( \eta \right)$, and $L_{k, r} \left( \zeta \right)$ are univariate B-spline basis functions of orders $p$, $q$, and $r$, corresponding to knot vectors $\boldsymbol{\Xi}$, $\boldsymbol{\mathcal{H}}$, and $\boldsymbol{\mathcal{Z}}$, respectively. The $\boldsymbol{B}_{i, j, k}$'s form a control mesh that does not conform to the actual geometry. To ensure that the generated control mesh corresponds to the B-spline volume, we adopt a mesh generation method \citep{otoguro2017space} based on the projection of a mesh generated with existing techniques to a B-spline volume.

For B-spline basis functions with $p=0$ and $p=1$, we obtain the same results as for standard piecewise constant and linear Lagrange basis functions, respectively. B-spline basis functions with $p \geq 2$, however, have several distinct features from their Lagrange-based counterparts.

\begin{figure}[t]
	\centering
	\includegraphics[bb=0 0 1027 346,width=12cm]{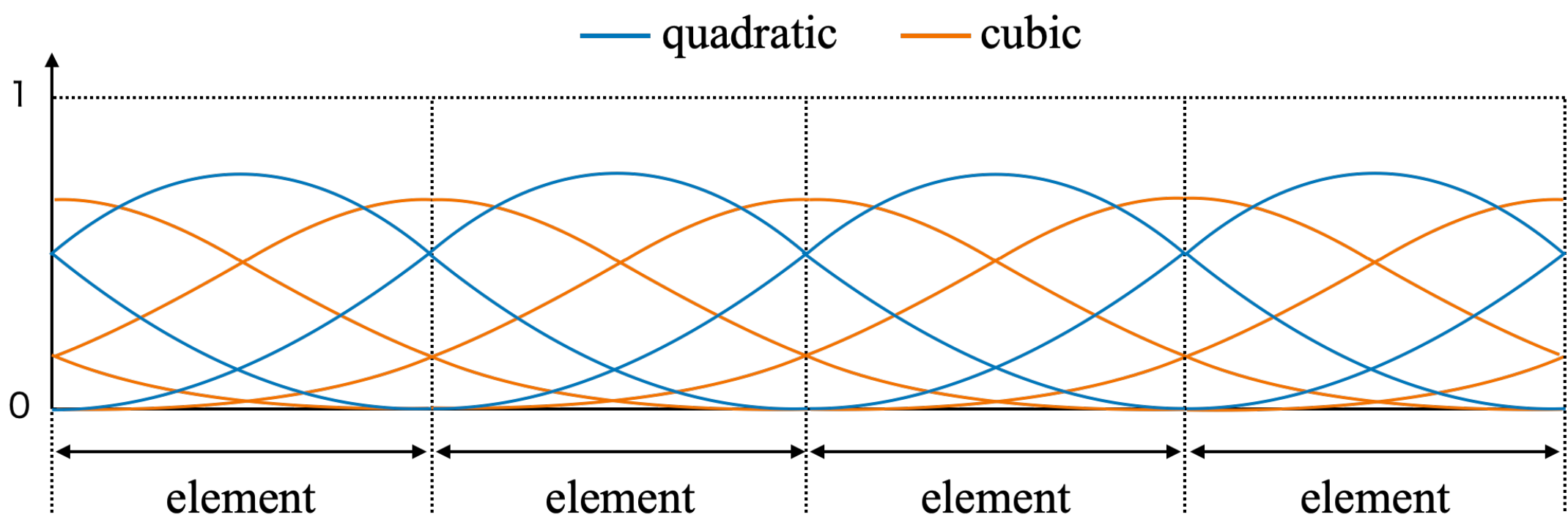}
	\caption{Quadratic and cubic B-spline basis functions for uniform knot vectors.}
	\label{fig:b-spline_function}
\end{figure}

The most important feature noted in the present work is that B-spline basis functions have higher continuity across the element boundaries than Lagrange basis functions, as shown in Figure \ref{fig:b-spline_function}. 
The figure depicts quadratic and cubic B-spline basis functions for uniform knot vectors, which are assembled by knots that are equally-spaced in the parametric space.
In general, $p$th order B-spline basis functions have $p-m_i$ continuous derivatives at knot $\xi_i$, where $m_i$ is the multiplicity of $\xi_i$ in the knot vector. In the present work, the multiplicity of all interior knots is defined as $m_i = 1$, where $p+2 \leq i \leq n$.
Thus, $p$th order B-spline basis functions have $C^{p-1}$-continuity across the interior element boundaries.

For numerical integration in B-spline meshes, Gaussian quadrature can be employed without any additional technique and is defined on each individual knot span in the parametric space.

A knot vector is said to be open if its first and last knot values appear $p+1$ times. 
B-spline basis functions formed from open knot vectors are interpolatory at the ends of the parameter space interval $\lbrack \xi_1, \xi_{n+p+1} \rbrack$. In general, B-spline basis functions are not interpolatory at interior knots. The present work employs open knot vectors.

\subsection{Proposed method: B-spline based s-version of finite element method (BSFEM)}
\label{sec:details_bsfem}
As shown in Section \ref{sec:difficulties_in_lagrange}, the conventional SFEM has problems with accuracy and computation time of numerical integration, as well as matrix singularity resulting in poor convergence for solving linear equations.

The former occurs because the global basis functions have low continuity across element boundaries. In conventional SFEM, Lagrange functions are used as both global and local basis functions. Because these functions have $C^0$-continuity, their derivatives exhibit discontinuities across element boundaries. Thus, the integrands often become discontinuous and the accuracy of Gaussian quadrature degrades when a local element contains several global elements.
The latter problem arises because the global and local basis functions are not guaranteed to be linearly independent of each other. Although many existing studies addressed this problem by constraining the degrees of freedom, the approaches employed therein are ad-hoc or incur additional computation time.

To solve the problem of accuracy, this study employs basis functions with high continuity across element boundaries as the global basis functions. Specifically, we define $p$th order B-spline basis functions ($p \geq 3$) as the global basis functions. The resulting integrands, including the first derivatives of basis functions, are smooth and continuous. Thus, Gaussian quadrature can be applied to the submatrices $\boldsymbol{K}^{\mathrm{GL}}$ and $\boldsymbol{K}^{\mathrm{LG}}$ without requiring additional techniques, ensuring accurate and efficient integration.
To solve the latter problem, we apply different functions to the global and local basis functions to guarantee their mutual linear independence. Specifically, we employ B-spline and Lagrange functions as the global and local basis functions, respectively. This approach is a potentially versatile and fundamental method to guarantee the uniqueness of the numerical solution.

\begin{figure}[t]
	\centering
	\includegraphics[bb=0 0 1279 763,width=10cm]{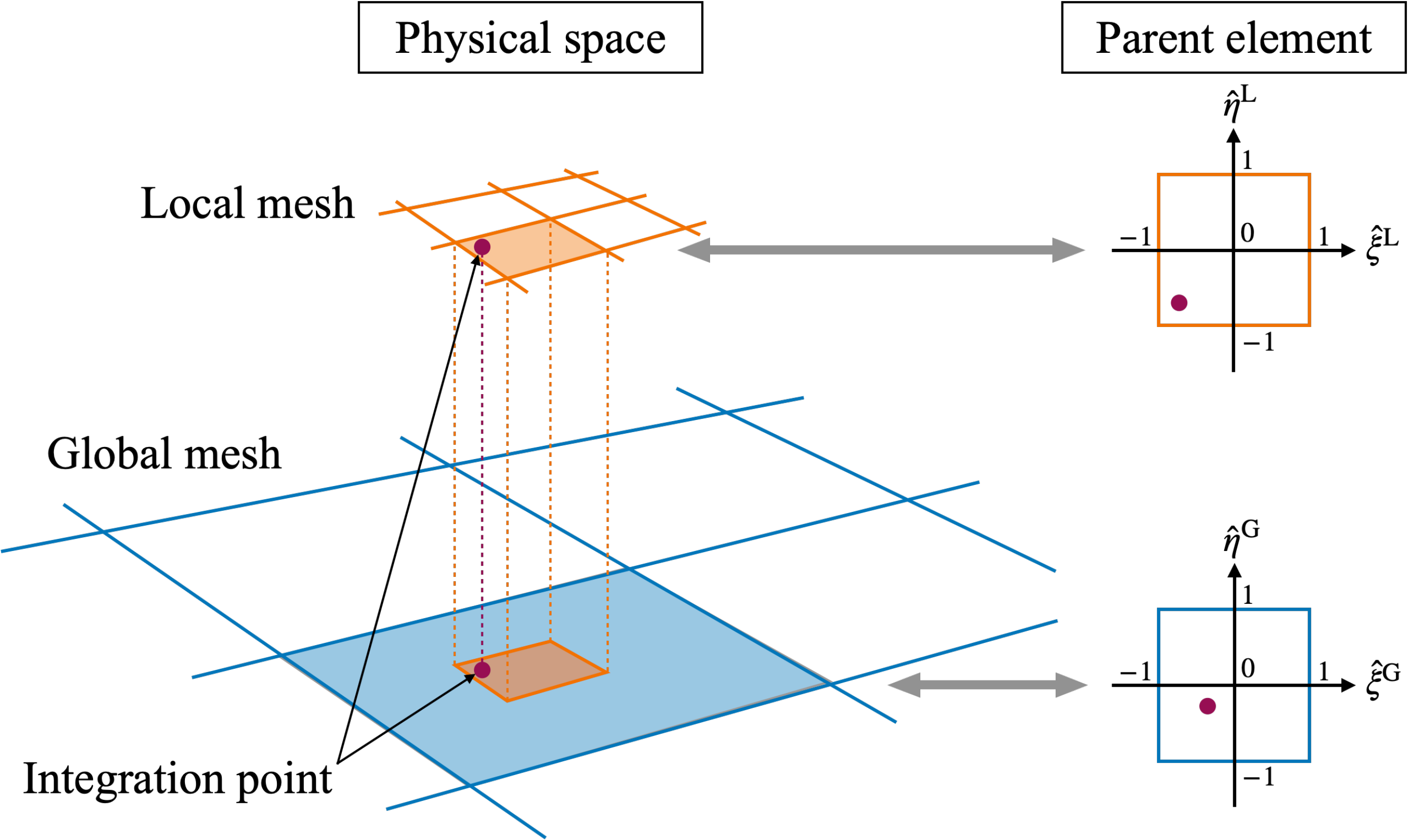}
	\caption{Schematic illustration of mapping between two meshes in SFEM.}
	\label{fig:sfem_mapping}
\end{figure}

As shown in Figure \ref{fig:sfem_mapping}, when the integrands are evaluated on the quadrature points defined in the parent element corresponding to the local element for the calculation of $\boldsymbol{K}^{\mathrm{GL}}$ and $\boldsymbol{K}^{\mathrm{LG}}$, the locations of these points in the parent element corresponding to the global element must be identified.
Because the problem of finding corresponding locations in the global mesh is nonlinear, an iterative procedure is generally required. When the shapes of the global elements are generated irregularly and the number of local elements increases, the mapping calculation may be rather inefficient and time-consuming.
By contrast, we employ a structured mesh as the B-spline based global mesh, assuming the imposition of boundary conditions by interface capturing approaches \citep{peskin1972flow, wang2004extended, zhang2004immersed, glowinski1999distributed, wagner2001extended, terada2003finite, udaykumar1996elafint, dunne2006eulerian}.
Thereby, the mapping in our proposed framework can be explicitly obtained without using iterative procedures, allowing us to improve computational efficiency.

\section{Verification of the proposed method}
\label{sec:verification}

\subsection{Target problem}
\label{sec:target_problem}
As described in Section \ref{sec:bsfem}, the conventional SFEM faces challenges in terms of numerical integration and independency of basis functions. To address these issues, we propose B-spline based SFEM, wherein B-spline and Lagrange basis functions are applied as the global and local basis functions, respectively.
The proposed method employs $p$th order B-spline basis functions ($p=2, 3$) as global basis functions, and $q$th order Lagrange basis functions ($q=1, 2, 3$) as local basis functions. That is, six pairs of global and local basis functions are tested.
In contrast, the conventional method employs $p$th order Lagrange basis functions ($p=1, 2, 3$) as global basis functions, and $q$th order Lagrange basis functions ($q=1, 2, 3$) as local basis functions. That is, nine pairs of global and local basis functions are tested.

\begin{figure}[t]
	\centering
	\includegraphics[bb=0 0 1671 728, width=12cm]{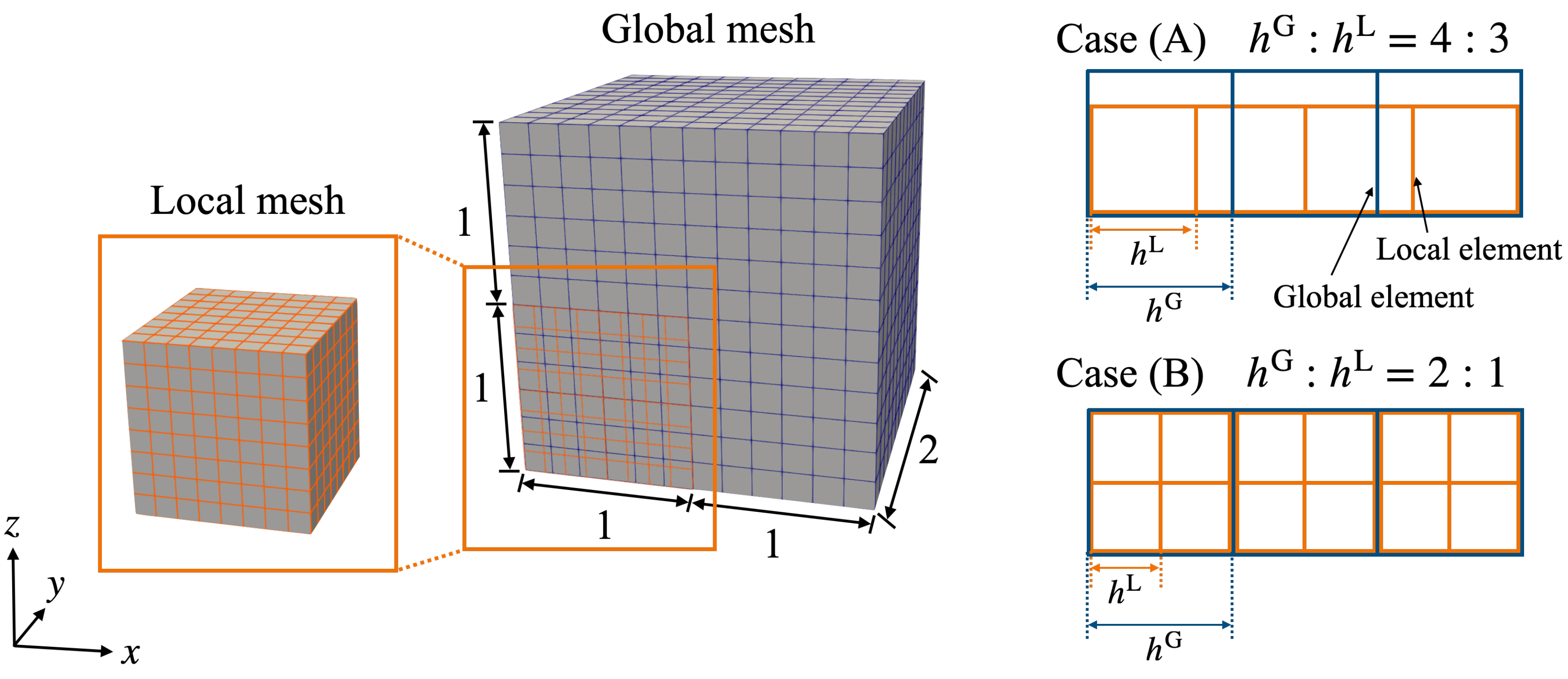}
	\caption{An example mesh model for SFEM.}
	\label{fig:analysis_mesh}
\end{figure}

The target problem is the basic Poisson's equation with Dirichlet boundary conditions as expressed in Section \ref{sec:formulation_of_s-fem}.
The domain for the analysis is defined as $[0, 2]^3$.
The global and local meshes with their boundary conditions are located in $[0, 2]^3$ and $[0, 1]^3$, respectively, considering each element as a cube.
An example of global and local meshes is illustrated in Figure \ref{fig:analysis_mesh}.

For the quantitative comparison, the same meshes were employed for the proposed and conventional methods in all verification tests.
The element sizes $h^{\mathrm{G}}$ and $h^{\mathrm{L}}$ in the respective meshes were used as evaluation parameters. 
Each test was performed under two conditions:
(A) $h^{\mathrm{G}}:h^{\mathrm{L}}=4:3$ and 
(B) $h^{\mathrm{G}}:h^{\mathrm{L}}=2:1$, as shown in Figure \ref{fig:analysis_mesh}.

The $n$-point Gaussian quadrature rule ensures the exact integration of a polynomial of degree $2n-1$ or lower. According to this rule, in FEM, using $p$th order basis functions, $p+1$-point Gaussian quadrature allows sufficiently accurate numerical integration.
In SFEM, if polynomials of degree $p$ are applied to the basis functions of one mesh and those of degree $q$ are applied to the basis functions of the other mesh $(p \geq q)$, all numerical integrals can be solved accurately using the $p+1$-point Gaussian quadrature. In Case (A), however, the global element boundary is contained within the local element, and a discontinuous function is integrated. Therefore, we employ high-order Gaussian quadrature to avoid the integral discontinuities that occur in the conventional SFEM \citep{fish1992s,lee2004combined,nakasumi2008crack,kishi2020dynamic,sawada2010high}.
In this study, we use the $p+8$-point Gaussian quadrature for Case (A) if polynomials of degree $p$ and $q$ are applied to the basis functions of two meshes, respectively $(p \geq q)$.
We evaluated the change in the relative $L^2$ error norm when the order of the Gaussian quadrature varied by 1 from $p+1$-point to $p+10$-point for each pair of basis functions. These tests were conducted for cases with the coarsest mesh $(h^{\mathrm{G}}=0.166667)$, where the relative $L^2$ error norm was expected to be maximal.
As a result, for all pairs of basis functions, the change in relative $L^2$ error norm was less than 5\% of the overall value of the relative $L^2$ error norm when using $p+8$-point Gaussian quadrature, compared to when using $p+7$-point and $p+9$-point Gaussian quadrature. In other words, when the $p+8$-point Gaussian quadrature is used, the effect of the relative $L^2$ error due to the integration of discontinuous functions on the overall relative $L^2$ error norm is sufficiently small, and the result is considered sufficiently accurate. We employed the same Gaussian quadrature for both the proposed and conventional methods to ensure a fair comparison.

In Case (B), no global element boundary is contained within the local element and no integration of the discontinuous function occurs. 
Accordingly, we employ the $p+1$-point Gaussian quadrature without any additional techniques to improve accuracy.
On the other hand, all global element boundaries in the local domain overlap with local element boundaries; thus, the basis functions of both meshes are mutually dependent in the conventional SFEM \citep{ooya2009linear}.
Here, we don't use any additional approach to avoid matrix singularity for both the proposed and conventional methods.

As described in Section \ref{sec:formulation_of_b-spline}, we define B-spline based meshes as structured meshes wherein every element is a cube. Open knot vectors are employed, and the multiplicity of all interior knots is defined as $1$. 
Hence, $p$th order B-spline basis functions have $C^{p-1}$-continuity across the element boundaries.
To form control meshes corresponding to the B-spline volumes, we adopt the mesh generation method \citep{otoguro2017space}.

In the present work, we employed 
(1) the relative $L^2$ error norm, 
(2) the iterative method for solving linear equations, with the number of iterations, and 
(3) the positive definiteness of the matrix
as the verification parameter.

The manufactured solution approach \citep{roache1998verification} is employed so that a convergence study on the relative $L^2$ error norm can be performed.
The exact solution is given as
\begin{equation}
    u = \sin{2 \pi x}\sin{2 \pi y}\sin{2 \pi z} + 10 \quad \mathrm{in} \; \Omega,
    \label{eq:manufactured_solution}
\end{equation}
and the resulting Poisson's equation for verification is defined as
\begin{equation}
    \Delta{u} + 12 {\pi}^2 \sin{2 \pi x}\sin{2 \pi y}\sin{2 \pi z} = 0 \quad \mathrm{in} \; \Omega.
\end{equation}
By applying Eq. \eqref{eq:manufactured_solution} to the problem domain boundary $\Gamma_D$ as the Dirichlet boundary condition, Eq. \eqref{eq:manufactured_solution} can be regarded as the exact solution of the problem.
The relative $L^2$ error norm $\varepsilon_{L^2}$ in the analysis based on SFEM is expressed as
\begin{align}
    \varepsilon_{L^2} 
    &= \frac{
    \sqrt{\int_{\Omega^{\mathrm{G}}} \left| u^h - u \right|^2 d \Omega}
    }{
    \sqrt{\int_{\Omega^{\mathrm{G}}} \left| u \right|^2 d \Omega}
    } \nonumber \\
    &= \frac{
    \sqrt{
    \int_{\Omega^{\mathrm{G}}\backslash\Omega^{\mathrm{L}}} \left| \left( u^{\mathrm{G}} \right)^h - u \right|^2 d \Omega
    +
    \int_{\Omega^{\mathrm{L}}} \left| \left( u^{\mathrm{G}} \right)^h + \left( u^{\mathrm{L}} \right)^h - u \right|^2 d \Omega 
    }}{\sqrt{
    \int_{\Omega^{\mathrm{G}}} \left| u \right|^2 d \Omega 
    }},
    \label{eq:sfem_error}
\end{align}
where $u$ is the exact solution and $\left( u^{\mathrm{G}} \right)^h$ and $\left( u^{\mathrm{L}} \right)^h$ are the calculated solutions in global and local meshes, respectively.
To simplify the integral calculations, both meshes are located so that the edges of the local domain are aligned along the global element boundaries, as shown in Figure \ref{fig:analysis_mesh}.
Integrations over $\Omega^{\mathrm{G}}$ and $\Omega^{\mathrm{G}}\backslash\Omega^{\mathrm{L}}$ are calculated based on the global mesh, whereas that over $\Omega^{\mathrm{L}}$ is calculated based on the local mesh.

Krylov-subspace methods, such as the conjugate gradient method, are often used to solve linear systems. However, the conjugate gradient method can only be used to solve a symmetric positive definite matrix.
In the conventional SFEM, the matrix \eqref{eq:submatrix_K} has been considered to be symmetric and positive definite in previous studies \citep{Okada2004-co}. However, no strict verification test or detailed discussion on the definiteness of the matrix in SFEM have been performed. 
In this study, we assessed the positive definiteness of the matrix in the proposed and conventional methods.
One method for testing whether a symmetric matrix is positive definite is Cholesky factorization \citep{higham2009cholesky, zhan1996computing}.
As described by \citet{higham2009cholesky}, upon running the Cholesky factorization algorithm, a matrix is considered positive definite if the algorithm completes without encountering any negative or zero pivots, and not positive definite otherwise.
In other words, if Cholesky factorization succeeds, all of the eigenvalues are positive, and the matrix is positive definite and regular. On the other hand, if Cholesky factorization fails, the matrix has at least one non-positive eigenvalue and is not positive definite. When the matrix has 0 eigenvalues, its determinant becomes 0 and the matrix is singular. When the matrix has negative eigenvalues, its determinant is not guaranteed to be 0 and the matrix is also not guaranteed to be singular. Loss of positive definiteness and matrix singularity are separate outcomes.

We verified the number of iterations required to solve the linear equations because, in many large-scale and realistic problems, the time required to solve said equations accounts for a large portion of the overall computation time. In other words, reducing the number of iterations in solving linear equations significantly reduces the overall computation time.

In the present work, we used the general-purpose linear equation solver library "Monolithic non-overlapping / overlapping DDM based linear equation solver (monolis)" \citep{monolis}.
Table \ref{tab:sfem_solver} lists other conditions required to solve linear equations.
\begin{table}[t]
    \begin{center}
        \caption{Analytical conditions for matrix calculation.}
        \begin{tabular}{lr} \hline
            Linear solver & conjugate gradient (CG) method\\
            Preconditioner & diagonal scaling method\\
            Convergence criterion & $1.0 \times 10^{-10}$ \\
            Maximum number of iterations & degrees of freedom \\ \hline
        \end{tabular}
        \label{tab:sfem_solver}
    \end{center}
\end{table}
The conjugate gradient method terminates in at most $n$ iterations, where $n$ corresponds to the degrees of freedom in the matrix, if no rounding errors are encountered \citep{hestenes1952methods}.
If the method fails to converge in $n$ iterations, we conclude that the matrix is singular. To minimize the influence of rounding errors, we employed the diagonal scaling method as a preconditioner.
Incidentally, $n$ is the sum of degrees of freedom of the global and local meshes.

In this paper, we tested the positive definiteness of the matrix via Cholesky factorization, and matrix singularity via the convergence of the conjugate gradient method.

\subsection{Results}
\label{sec:results}
\begin{figure}[t]
	\centering
	\includegraphics[bb=0 0 660.38 464.82, width=14cm]{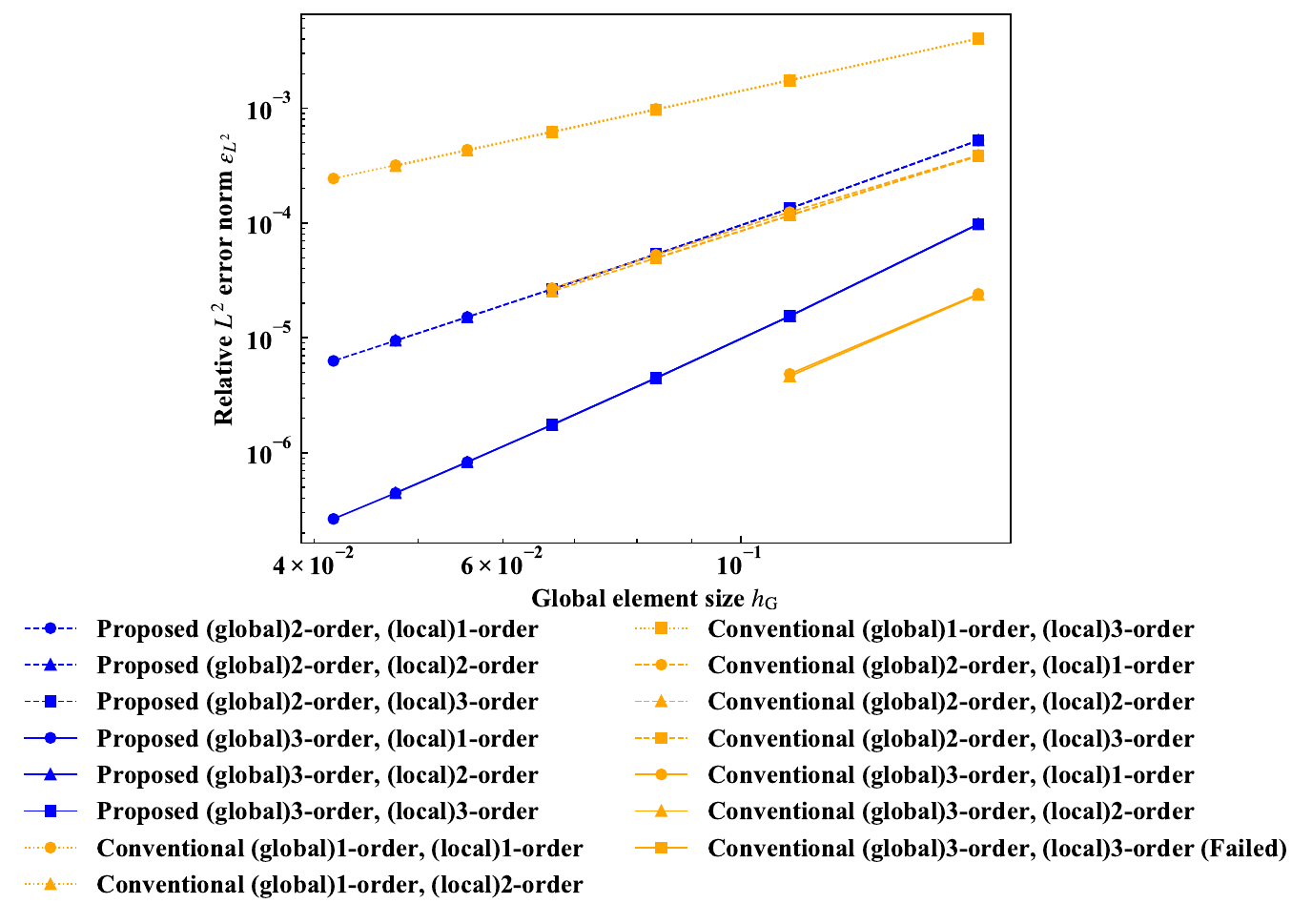}
	\caption{Convergence of relative $L^2$ error norm against global element size $h^{\mathrm{G}}$ for Case (A). "(Failed)" indicates that solutions did not converge in all cases except where calculation was not possible due to insufficient memory.}
	\label{fig:L2error_convergence_caseA}
\end{figure}

\begin{figure}[t]
	\centering
	\includegraphics[bb=0 0 706.94 464.82, width=14cm]{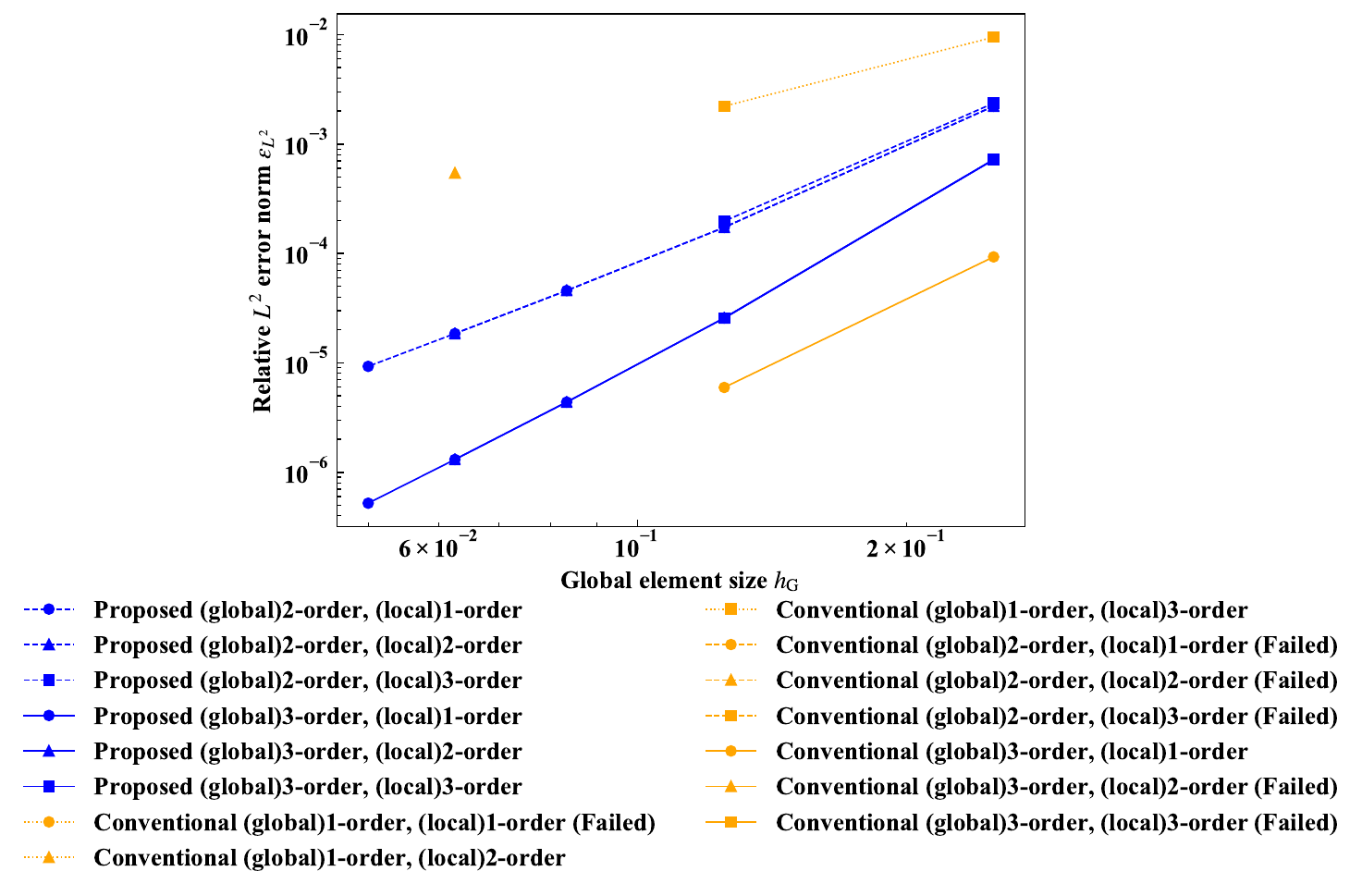}
	\caption{Convergence of relative $L^2$ error norm against global element size $h^{\mathrm{G}}$ for Case (B). "(Failed)" indicates that solutions did not converge in all cases except where calculation was not possible due to insufficient memory.}
	\label{fig:L2error_convergence_caseB}
\end{figure}

\begin{figure}[t]
	\centering
	\includegraphics[bb=0 0 660.38 464.82, width=14cm]{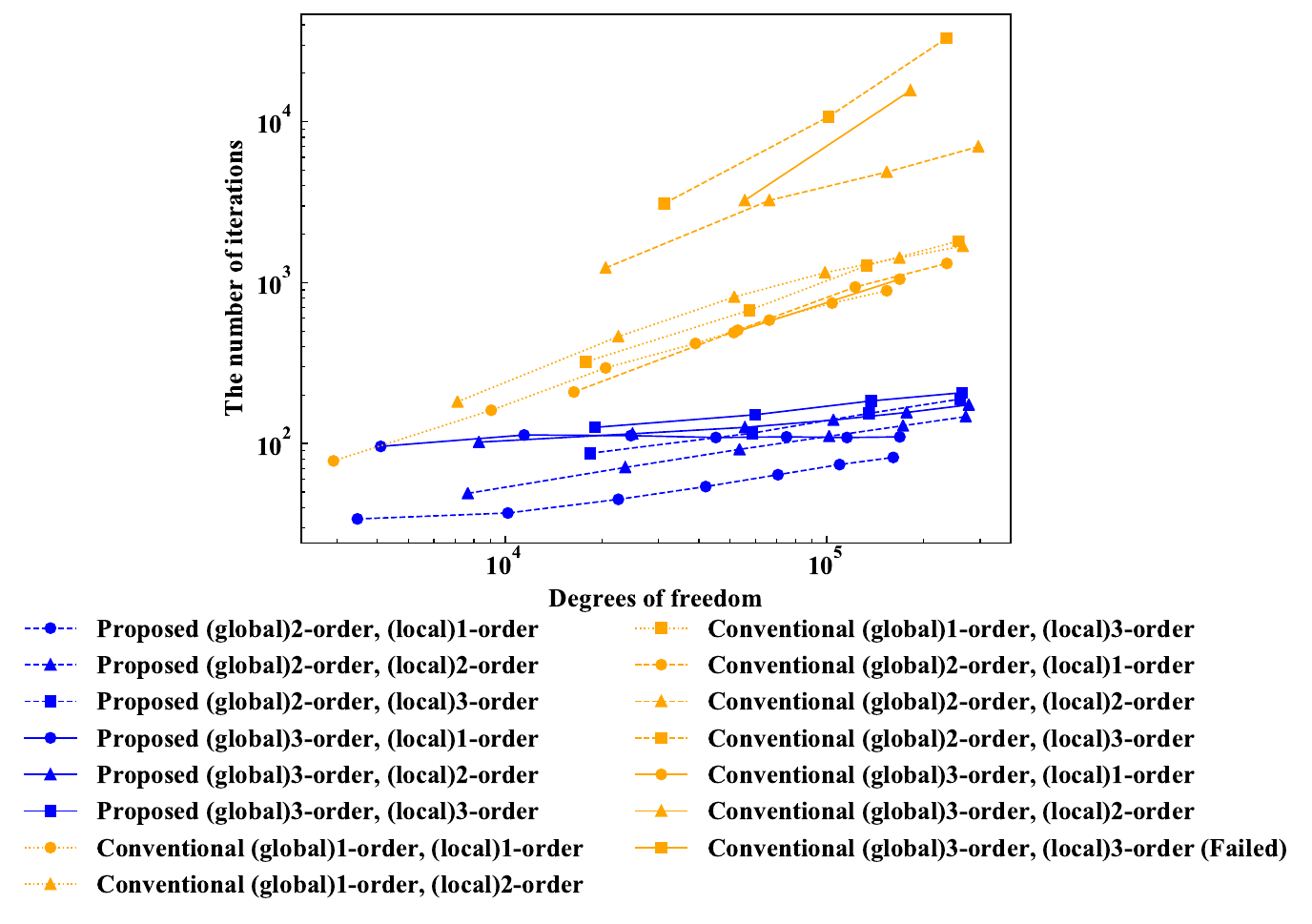}
	\caption{Number of iterations until convergence against degrees of freedom for Case (A). "(Failed)" indicates that solutions did not converge in all cases except where calculation was not possible due to insufficient memory.}
	\label{fig:iteration_dof_caseA}
\end{figure}

\begin{figure}[t]
	\centering
	\includegraphics[bb=0 0 706.94 464.82, width=14cm]{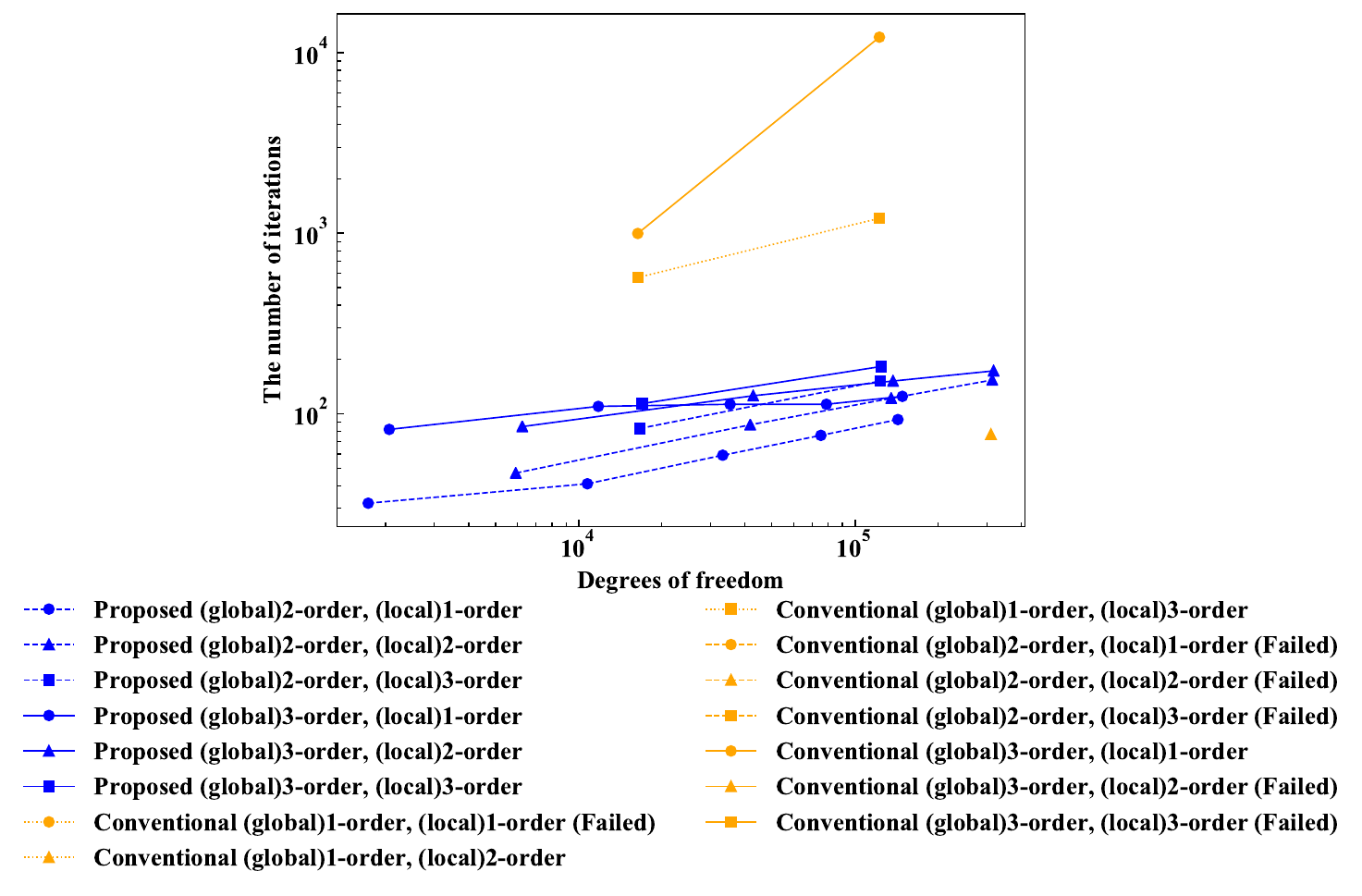}
	\caption{Number of iterations until convergence against degrees of freedom for Case (B). "(Failed)" indicates that solutions did not converge in all cases except where calculation was not possible due to insufficient memory.}
	\label{fig:iteration_dof_caseB}
\end{figure}

\begin{figure}[t]
	\centering
	\includegraphics[bb=0 0 660.38 464.82, width=14cm]{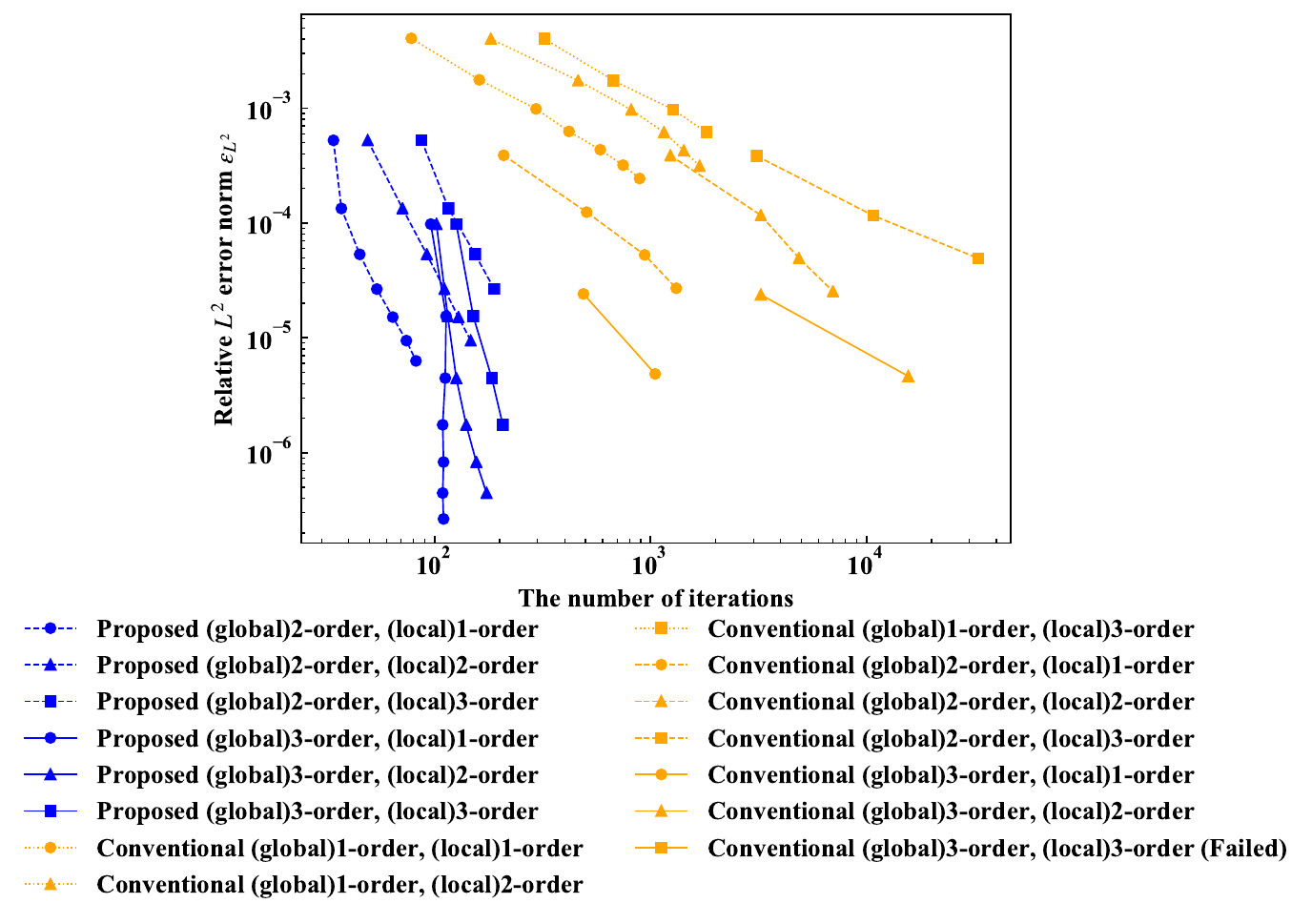}
	\caption{Relative $L^2$ error norm against number of iterations until convergence for Case (A). "(Failed)" indicates that solutions did not converge in all of its cases except cases where calculation was not possible due to insufficient memory.}
	\label{fig:L2error_iteration_caseA}
\end{figure}

\begin{figure}[t]
	\centering
	\includegraphics[bb=0 0 706.94 464.82, width=14cm]{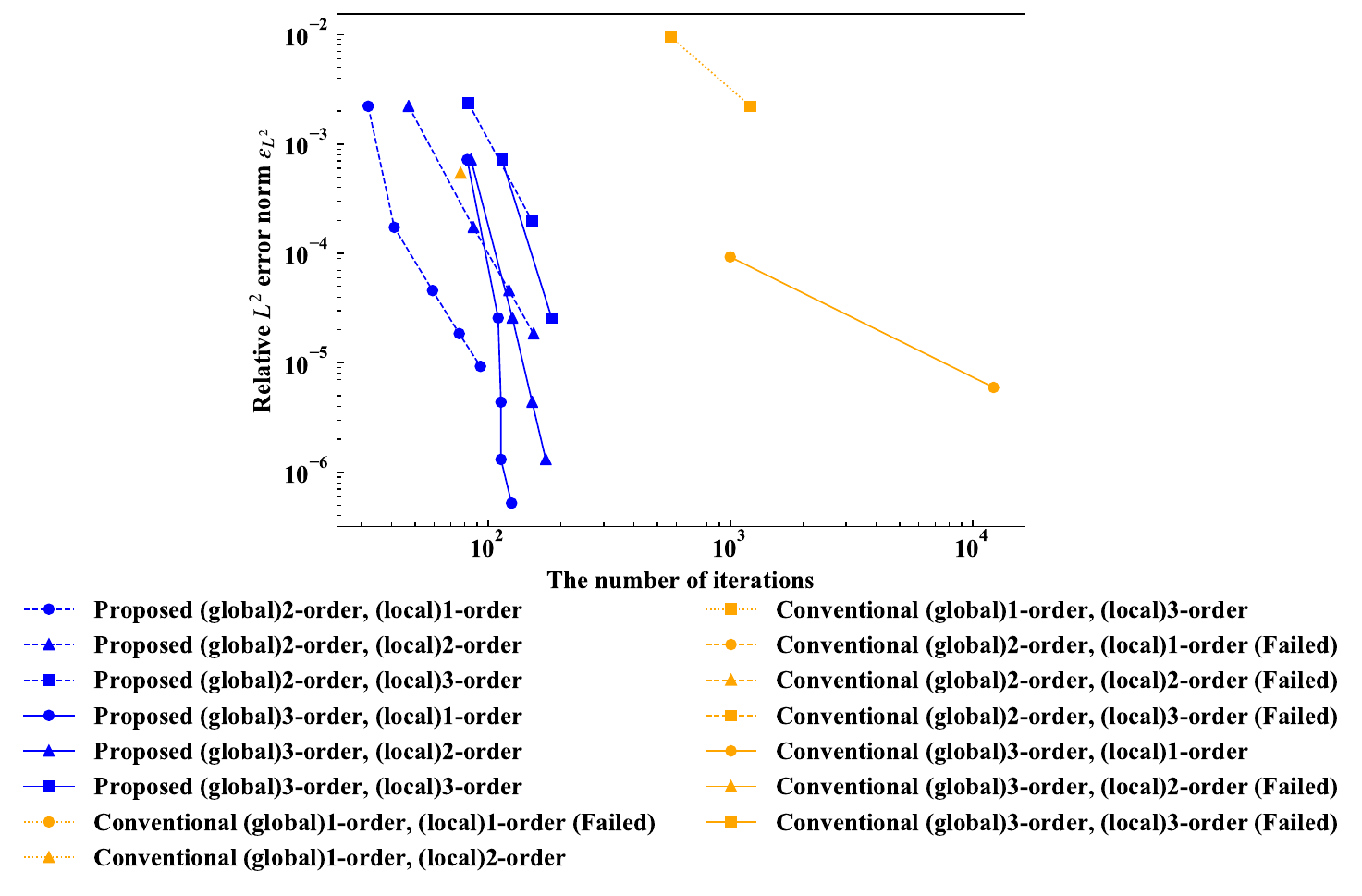}
	\caption{Relative $L^2$ error norm against number of iterations until convergence for Case (B). "(Failed)" indicates that solutions did not converge in all cases except where calculation was not possible due to insufficient memory.}
	\label{fig:L2error_iteration_caseB}
\end{figure}

In this verification, we discuss the differences in the overall trends of the results of the proposed and conventional methods with respect to the relative $L^2$ error norm, convergence of the conjugate gradient method, and positive definiteness of the matrix. In the following graphs 
(Figures \ref{fig:L2error_convergence_caseA}, \ref{fig:L2error_convergence_caseB}, \ref{fig:iteration_dof_caseA}, \ref{fig:iteration_dof_caseB}, \ref{fig:L2error_iteration_caseA}, and \ref{fig:L2error_iteration_caseB}),
the results of the proposed method are shown in blue and those of the conventional method in orange.

The convergence of the relative $L^2$ error norm defined in Eq. \eqref{eq:sfem_error} against the global element size $h^{\mathrm{G}}$ was evaluated in the proposed and conventional methods. The numerical results for Cases (A) and (B) are shown in Figures \ref{fig:L2error_convergence_caseA} and \ref{fig:L2error_convergence_caseB}, respectively.
The results show that the proposed method exhibits better error convergence in all cases when compared with the same order basis functions in the conventional method.
Thus, the proposed method seems to show better accuracy than the conventional method for larger-scale analysis that requires more detailed meshes.
Comparing the results of Cases (A) and (B) for each method, the errors are comparable in all cases. This indicates that the high-order Gaussian quadrature method can be used to calculate Case (A) of the conventional method, where the integration of discontinuous functions occurs, with sufficient accuracy. However, this approach incurs significant computation time, in line with other methods of improving accuracy which are more complex.
In addition, both methods are most accurate when the global basis functions are third order, and there is little variability depending on the local basis functions. This seems to be a natural result because the global basis functions have a dominant impact on the relative $L^2$ error distribution over the entire domain, and this problem does not generate local regions that require high resolution.
We note that when using the conventional method in most cases of (B), solutions did not converge even after the maximum number of iterations was reached, as shown in Table \ref{tab:sfem_solver}.
These results are seemingly caused by matrix singularity.
By contrast, using the proposed method, solutions converged in all cases of (A) and (B). These results indicate that the basis functions in the proposed method are guaranteed to be mutually linearly independent, thereby avoiding matrix singularity.

\begin{figure}[t]
  \begin{minipage}[b]{0.35\linewidth}
    \centering
    \includegraphics[bb=0.000000 0.000000 832.000000 672.000000,scale=0.17]{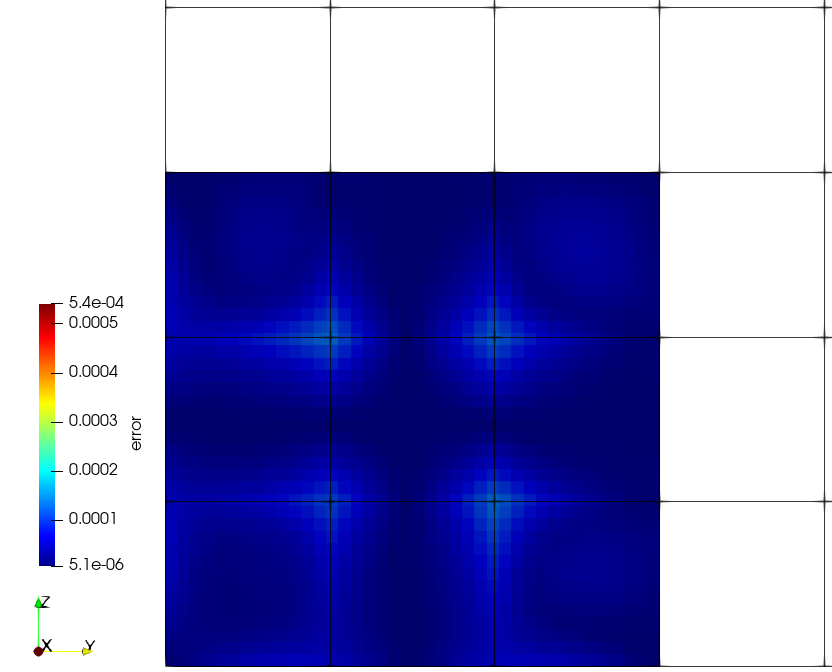}
  \end{minipage}
  \hspace{0.01\columnwidth}
  \begin{minipage}[b]{0.3\linewidth}
    \centering
    \includegraphics[bb=0.000000 0.000000 676.000000 674.000000,scale=0.17]{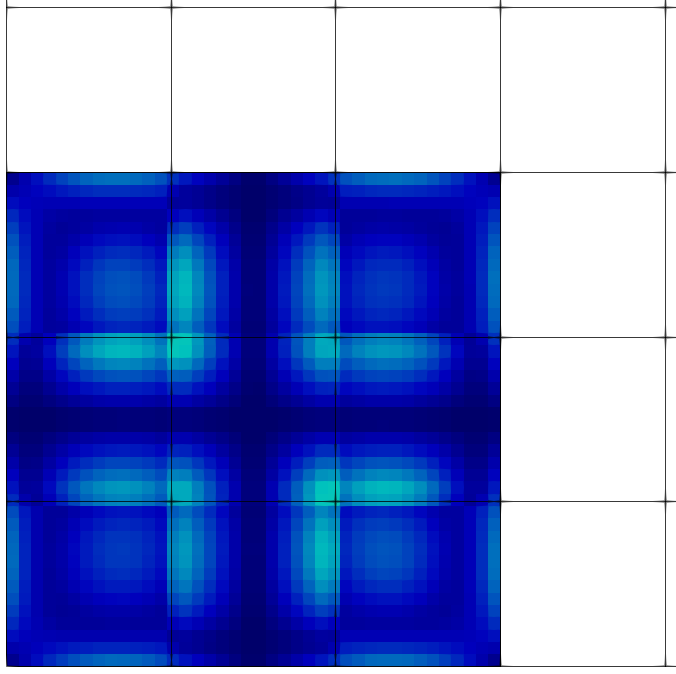}
  \end{minipage}
  \hspace{0.01\columnwidth}
  \begin{minipage}[b]{0.3\linewidth}
    \centering
    \includegraphics[bb=0.000000 0.000000 676.000000 670.000000,scale=0.17]{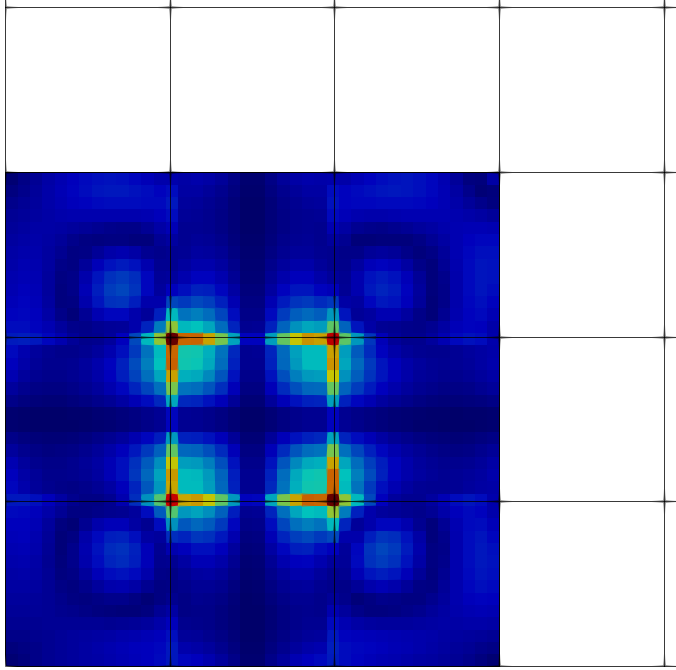}
  \end{minipage}
  \caption{Cross-sectional view (yz plane, x=0.33) of the relative $L^2$ error distribution in the local domain for the case of $h_\mathrm{G}:h_\mathrm{L}=40:3$. Global basis functions are (left:) cubic B-spline basis functions, (middle:) quadratic B-spline basis functions, and (right:) cubic Lagrange basis functions. The contours of error are shown in color (min: $1.2 \times 10^{-7}$, max: $5.4\times 10^{-4}$).}
  \label{fig:error_distribution}
\end{figure}

To qualitatively assess the error due to the integration of discontinuous functions, the relative $L^2$ error distribution in the local domain was verified for three cases with different continuity of the global basis functions, where the ratio of global to local element sizes was extreme: $h_{\mathrm{G}}:h_{\mathrm{L}}= 40:3$.
This test was performed under the following three cases: (1) cubic B-spline basis functions with $C^2$-continuity across element boundaries, (2) quadratic B-spline basis functions with $C^1$-continuity across element boundaries, and (3) cubic Lagrange basis functions with $C^0$-continuity across element boundaries are applied as global basis functions.
In all cases, linear Lagrange functions were employed as local basis functions.
To evaluate the error due to the integration of discontinuous functions, we applied the $p+1$-point Gaussian quadrature without any additional techniques to improve accuracy when $p$th and $q$th order basis functions were applied to two meshes, respectively $(p \geq q)$.
The results are shown in Figure \ref{fig:error_distribution}.
The distribution of the relative $L^2$ error norm for each element in the local domain is shown by the colored contour, and the element boundaries of the global mesh are denoted by black lines.
These results indicate that the case wherein the global basis functions have lower continuity produces significantly larger errors for local elements that are located across global element boundaries, and that the error distribution is discontinuous in proximity of said boundaries.
In contrast, the proposed method, which uses cubic B-spline functions for the global basis, exhibits almost no such errors. 
That is, the proposed method allows sufficiently accurate computation using Gaussian quadrature without any additional and computationally expensive techniques to improve accuracy.
Consequently, the proposed method can further reduce computation time at the same level of accuracy.

Relationships between the number of iterations required for convergence and the degrees of freedom for Cases (A) and (B) are shown in Figures \ref{fig:iteration_dof_caseA} and \ref{fig:iteration_dof_caseB}, respectively.
The results show that using the proposed method, the solutions converged at a small number of iterations in all cases of (A) and (B) even with large degrees of freedom. On the other hand, using the conventional method, the number of iterations increased significantly with degrees of freedom for both cases. Furthermore, many cases of (B) failed to converge using the conventional method, indicating matrix singularity.
Therefore, not only does the conventional method fail to converge in some cases, but convergence is also slow and computationally intensive when solving simultaneous linear equations in almost all cases. On the other hand, the proposed method exhibits excellent convergence in all cases of (A) and (B). 
From these results, we conclude that the proposed method guarantees linear independence of basis functions, and thus, it has excellent convergence.
Reducing the number of iterations to solve linear equations implies reducing the overall computation time. In addition, the proposed method does not require computationally expensive or ad-hoc techniques to improve convergence \citep{fish1992s,FISH1994135,ANGIONI2011780,ANGIONI2012559,yue2005adaptive,yue2007adaptive,fan2008rs,nakasumi2008crack,park2003efficient,ooya2009linear}, further reducing computation time and simplifying the meshing procedure of SFEM.

Figures \ref{fig:iteration_dof_caseA} and \ref{fig:iteration_dof_caseB} show that the conventional method converged very slowly, even in cases where the solution converged and the matrix was not considered singular. For further verification of this problem, we focused on Case A and tested the positive definiteness of the matrices.
Table \ref{tab:definiteness_caseA} lists the verification results of the positive definiteness of the matrices in the proposed and conventional SFEM for Case (A).
We conclude that matrices that can be Cholesky decomposed are positive definite matrices, whereas those that cannot be Cholesky decomposed are not positive definite matrices.
In the tables, "Pass" means that decomposition succeeded and "Fail" means that it failed.
Blank columns indicate cases where calculation was not possible due to insufficient memory, which are discussed further.
\begin{table}[t]
    \begin{center}
        \caption{Positive definiteness of the matrices in proposed and conventional methods for Case (A).}
        \scalebox{0.8}{
            \begin{tabular}{cccccccccc} \hline
                & Order of & Order of & \multicolumn{7}{c}{Number of global elements}\\
                & global basis & local basis & $12^3$ & $18^3$ & $24^3$ & $30^3$ & $36^3$ & $42^3$ & $48^3$\\ \hline
                \multirow{6}{*}{
                    \begin{tabular}{c} Proposed \\ method\end{tabular}
                }
                & 2 & 1 & Pass & Pass & Pass & Pass & Pass & &\\ 
                & 2 & 2 & Pass & Pass & Pass &      &      & &\\
                & 2 & 3 & Pass & Pass &      &      &      & &\\
                & 3 & 1 & Pass & Pass & Pass & Pass & Pass & &\\ 
                & 3 & 2 & Pass & Pass & Pass &      &      & &\\
                & 3 & 3 & Pass & Pass &      &      &      & &\\\hline
                \multirow{9}{*}{
                    \begin{tabular}{c} Conventional \\ method\end{tabular}
                }
                & 1 & 1 & Fail & Fail & Fail & Fail & Fail & &\\ 
                & 1 & 2 & Fail & Fail & Fail &      &      & &\\
                & 1 & 3 & Fail & Fail &      &      &      & &\\
                & 2 & 1 & Fail & Fail &      &      &      & &\\
                & 2 & 2 & Fail & Fail &      &      &      & &\\
                & 2 & 3 & Fail &      &      &      &      & &\\
                & 3 & 1 & Fail &      &      &      &      & &\\
                & 3 & 2 & Fail &      &      &      &      & &\\
                & 3 & 3 &      &      &      &      &      & &\\\hline
            \end{tabular}
        }
        \label{tab:definiteness_caseA}
    \end{center}
\end{table}
The results show that the matrices in the proposed method were positive definite in all cases, whereas those in the conventional method were not positive definite in all cases for Case (A).
As noted above, the conjugate gradient method and many other iterative methods can only be used on positive definite matrices.
Thus, in conventional SFEM, not only can the matrix become singular, but it can also lose its positive definiteness.
This result is a novel finding that contradicts the assumptions of existing studies, suggesting a possible cause of the poor convergence of the iterative approach when solving linear equations in the conventional method.

Relationships between the relative $L^2$ error norm and the number of iterations until convergence for Cases (A) and (B) are shown in Figures \ref{fig:L2error_iteration_caseA} and \ref{fig:L2error_iteration_caseB}.
These results indicate that the proposed method requires far fewer iterations than the conventional method for the same relative $L^2$  error norm.
Because the computation time for solving linear equations is a huge component of the overall computation time, the proposed method with excellent convergence results in much less computation time at the same level of accuracy.

These results demonstrate that the proposed method achieves sufficient accuracy of numerical integration and excellent convergence without requiring additional computationally expensive or ad-hoc techniques.
We therefore conclude that the proposed method is expected to achieve the same level of accuracy with less computation time than the conventional method, or vice versa.

\section{Conclusions and future works}
\label{sec:conclusions}
This study proposed a B-spline based SFEM that fundamentally solved the challenges of the conventional SFEM in numerical integration and matrix singularity.
There are two challenges with the conventional method. First, the inaccuracy of numerical integration in the term representing the interaction between the global and local meshes. This problem arises when the global basis functions have low continuity across the element boundaries.
Second, because linear independence of the global and local basis functions is not guaranteed, the matrices become singular or nearly singular, and the convergence of solving simultaneous linear equations deteriorates.
Thus, this study fundamentally solved these problems by applying cubic B-spline basis functions with $C^2$-continuity across element boundaries as global basis functions, while retaining Lagrangian basis functions as local basis functions.

We verified the proposed and conventional SFEM using the relative $L^2$ error norm, number of iterations for solving linear equations, and positive definiteness of the matrix as parameters.
Our results indicate that the proposed method can be computed with sufficient accuracy using Gaussian quadrature without requiring additional and computationally expensive techniques to improve accuracy. Furthermore, the conventional method was observed to require many iterations for convergence, and occasionally failed to converge even at the maximum number of iterations. In contrast, the proposed method exhibited convergence at significantly small numbers of iterations for the same problems. The relationship between excellent and poor convergence in the proposed and conventional methods and the positive definiteness of the matrix was identified. These results indicate that the proposed method guarantees linear independence of basis functions and has excellent convergence.
Therefore, we concluded that the proposed method has potential to reduce computation time while maintaining accuracy.

This study presents the first steps toward improving localized mesh refinement schemes in interface-capturing approaches for moving boundary problems. Future efforts will be focused on the introduction of ALE schemes for moving local mesh tracking to the boundary layers and surrounding area, coupling with IB methods and other approaches to handle boundary conditions at interfaces, and extending it to unsteady nonlinear problems.

\section*{CRediT authorship contribution statement}
\textbf{Nozomi Magome:} Data curation, Formal analysis, Investigation, Methodology, Software, Validation, Visualization, Writing – original draft.
\textbf{Naoki Morita:} Methodology, Software, Supervision, Writing – review \& editing.
\textbf{Shigeki Kaneko:} Methodology, Supervision, Writing – review \& editing.
\textbf{Naoto Mitsume:} Conceptualization, Funding acquisition, Project administration, Resources, Supervision, Writing – review \& editing.

\section*{Declaration of competing interest}
The authors declare that they have no known competing financial interests or personal relationships that could have appeared to influence the work reported in this paper.

\section*{Data availability}
Data will be made available on request.

\section*{Acknowledgements}
This work was supported by JST FOREST Grant Number JPMJFR215S and JSPS KAKENHI Grant Numbers 22H03601, 23H00475.




\bibliographystyle{elsarticle-num-names} 
\bibliography{bibliography_paperpile}


\end{document}